\documentclass{article}
\usepackage{geometry}
\usepackage{amssymb}
\usepackage{amsmath,cases}
\usepackage{amsthm}
\usepackage{enumerate}
\usepackage{xcolor}
\usepackage{graphicx}
\usepackage{textcomp}
\usepackage{mathrsfs}
\usepackage{bbm}
\usepackage{bm}
\usepackage{caption}
\usepackage{chemarrow}
\usepackage{stmaryrd}
\usepackage{url}
\usepackage{calc}
\usepackage{ifthen}
\usepackage{psfrag}
\usepackage[active]{srcltx}   
\usepackage[all,cmtip]{xy}
\usepackage{pdfsync}
\usepackage[refpage]{nomencl}  
\usepackage{yhmath}
\usepackage{framed}
\usepackage[colorlinks=true,bookmarksnumbered=true, pdfauthor={Jérémie Bettinelli}]{hyperref}									

\renewcommand{\emptyset}{\varnothing}

\newcommand{\8}{\infty}

\newcommand{\bs}{\backslash}
\newcommand{\eps}{\varepsilon}

\newcommand{\de}{\mathrel{\mathop:}\hspace*{-.6pt}=}

\newcommand{\ooo}{\bullet}

\newcommand{\si}{\sigma}

\newcommand{\ton}{\xrightarrow[n\to\infty]{}}

\newcommand{\tol}{\xrightarrow[n \to \infty]{(d)}}
\newcommand{\tolk}{\xrightarrow[k \to \infty]{(d)}}

\newcommand{\ee}{\mathbbm{e}}
\newcommand{\E}[1]{\mathbb{E}\left[ #1 \right]}

\newcommand{\N}{\mathbb{N}}
\newcommand{\Pb}{\mathbb{P}}

\newcommand{\R}{\mathbb{R}}

\newcommand{\Sgp}{\mathbf\Sigma_{g,p}^{\partial}}
\newcommand{\Sgpof}[2]{\mathbf\Sigma_{#1,#2}^{\partial}}
\newcommand{\Z}{\mathbb{Z}}
\newcommand{\Zp}{\Z_+}

\newcommand{\cB}{\mathcal{B}}
\newcommand{\C}{\mathcal{C}}

\newcommand{\EE}{\mathcal{E}}
\newcommand{\F}{\mathcal{F}}
\newcommand{\calG}{\mathcal{G}}

\newcommand{\K}{\mathcal{K}}
\newcommand{\lL}{\mathcal{L}}
\newcommand{\LMnsun}{\mathfrak{M}_{n,\sigma}^1}
\newcommand{\LMnsnun}{\mathfrak{M}_{n,\sigma_n}^1}
\newcommand{\LMnsde}{\mathfrak{M}_{n,\sigma}^2}
\newcommand{\LMnsnde}{\mathfrak{M}_{n,\sigma_n}^2}
\newcommand{\Mnsun}{\mathcal{M}_{n,\sigma}^1}
\newcommand{\Mnsde}{\mathcal{M}_{n,\sigma}^2}
\newcommand{\M}{\mathcal{M}}

\newcommand{\Q}{\mathcal{Q}}

\newcommand{\Qns}{\mathcal{Q}_{n,\sigma}}
\newcommand{\Qnsn}{\mathcal{Q}_{n,\sigma_n}}

\newcommand{\rR}{\mathcal{R}}

\newcommand{\X}{\mathcal{X}}

\newcommand{\NC}{\mathcal{N}}

\newcommand{\Mi}{\mathscr{M}}

\newcommand{\bb}{\mathfrak{b}}
\newcommand{\fB}{\mathfrak{B}}

\newcommand{\f}{\mathfrak{f}}
\newcommand{\fF}{\mathfrak{F}}

\newcommand{\fI}{\mathfrak{I}}
\newcommand{\lab}{\mathfrak{l}}
\newcommand{\Lab}{\mathfrak{L}}

\newcommand{\m}{\mathfrak{m}}
\newcommand{\mM}{\mathfrak{M}}
\newcommand{\q}{\mathfrak{q}}
\newcommand{\qi}{\mathfrak{q}_\infty}
\newcommand{\qis}{\mathfrak{q}_\infty^\sigma}

\newcommand{\s}{\mathfrak{s}}
\newcommand{\Sg}{\mathfrak{S}}
\newcommand{\tr}{\mathfrak{t}}

\newcommand{\cov}{\operatorname{cov}}

\newcommand{\un}[1]{\mathbbm{1}_{#1}}

\newcommand{\st}[2]{[#1 \rightarrowtriangle #2]}
\newcommand{\ori}{\check}

\newcommand{\bp}{\rho^\bullet}
\newcommand{\eqq}{\sim}
\newcommand{\eqt}{\simeq}

\newcommand{\fl}{f\hspace{-0.6mm}l}

\newcommand{\dish}{\delta_\mathcal{H}}

\newcommand{\dGH}{d_{\operatorname{GH}}}

\newcommand{\di}{d_\infty}
\newcommand{\disig}{d_\infty^\sigma}
\newcommand{\pii}{\pi}

\newcommand{\BDG}{Bouttier--Di~Francesco--Guitter }
\newcommand{\g}{(8n/9)^{1/4}}

\newcommand{\vh}{\rlap{$\vec{\hspace{2mm}}$}{h}}
\newcommand{\Geo}{\operatorname{Geod}}
\newcommand{\Mult}{\operatorname{Mult}}
\newcommand{\HMult}{\operatorname{HMult}}
\newcommand{\dpath}{d_{\operatorname{path}}}

\newcommand{\sand}{\qquad\text{ and }\qquad}

\newcommand{\lhb}{[[}

\newcommand{\rhb}{]]}

\newcommand{\lp}{\left(}
\newcommand{\rp}{\right)}
\newcommand{\lb}{\left\{}
\newcommand{\rb}{\right\}}

\newcommand{\lt}{\left|}
\newcommand{\rt}{\right|}

\definecolor{gris}{gray}{0.7}
\definecolor{grisf}{gray}{0.4}
\definecolor{vert}{rgb}{0,.5547,0}

\theoremstyle{plain}
\newtheorem{thm}{Theorem}
\newtheorem{lem}[thm]{Lemma}
\newtheorem{prop}[thm]{Proposition}
\newtheorem{conj}{Conjecture}
\newtheorem{corol}[thm]{Corollary}
\newtheorem{defi}{Definition}

\theoremstyle{definition}
\newtheorem*{rem}{Remark}
\newtheorem{rems}{Remark}

\newenvironment{pre}[1][\proofname]{%
  \proof[#1]%
}{\endproof}

\makenomenclature

\makeatletter
\newcommand \Dotfill {\leavevmode \leaders \hb@xt@ .72em{\hss .\hss } \hfill \kern \z@}
\makeatother
\renewcommand*{\pagedeclaration}[1]{\Dotfill\ifthenelse{#1<10}{\hspace{.72em}}{\hspace{.2em}}\hyperpage{#1}}
\makeatletter
\renewcommand*{\@fnsymbol}[1]{\ensuremath{\ifcase#1\or \dagger\or \ddagger\or
   \mathsection\or \mathparagraph\or \|\or **\or \dagger\dagger
   \or \ddagger\ddagger \else\@ctrerr\fi}}
\makeatother
\title{Geodesics in Brownian surfaces (Brownian maps)}
\author{J\'er\'emie \textsc{Bettinelli}\thanks{CNRS \& Institut \'Elie Cartan de Lorraine; \href{mailto:jeremie.bettinelli@normalesup.org}{\nolinkurl{jeremie.bettinelli@normalesup.org}}; \nolinkurl{www.normalesup.org/}\texttildelow\nolinkurl{bettinel}.}}

\begin{document}
\maketitle

\begin{abstract}
We define a class a metric spaces we call \emph{Brownian surfaces}, arising as the scaling limits of random maps on general orientable surfaces with a boundary and we study the geodesics from a uniformly chosen random point. These metric spaces generalize the well-known \emph{Brownian map} and our results generalize the properties shown by Le~Gall on geodesics in the latter space. We use a different approach based on two ingredients: we first study typical geodesics and then all geodesics by an ``entrapment'' strategy. In particular, we give geometrical characterizations of some subsets of interest, in terms of geodesics, boundary points and concatenations of geodesics forming a loop that is not homotopic to~$0$.

\bigskip
\begin{center}
\textbf{R\'esum\'e}
\end{center}

On d\'efinit une classe d'espaces m\'etriques al\'eatoires que nous appelons \emph{surfaces browniennes}~: ces objets apparaissent comme limites d'\'echelle de cartes al\'eatoires sur des surfaces orientables \`a bord g\'en\'erales. Dans un second temps, on \'etudie les g\'eod\'esiques \'emanant d'un point choisi uniform\'ement au hasard. Les surfaces browniennes g\'en\'eralisent la fameuse \emph{carte brownienne} et nos r\'esultats g\'en\'eralisent les propri\'et\'es obtenues par Le~Gall sur les g\'eod\'esiques dans cet espace. On utilise une approche diff\'erente reposant sur deux ingr\'edients~: on \'etudie d'abord les g\'eod\'esiques aux points typiques et on attrape ensuite les autres points en les ``encerclant'' par de telles g\'eod\'esiques. En particulier, on obtient des caract\'erisations g\'eom\'etriques de certains sous-ensembles d'int\'er\^et en termes de g\'eod\'esiques, points du bord et concat\'enations des g\'eod\'esiques formant une boucle non homotope \`a~$0$.

\paragraph{Key words and phrases:} Brownian surfaces; Brownian map; geodesics; random maps; scaling limits; Gromov--Hausdorff topology; random metric spaces; bijections.

\paragraph{AMS 2000 Subject Classification:} 60F17; 60C05; 60D05; 57N05; 05C80; 05C12.

\end{abstract}

\begin{figure}[ht!]
	\centering\href{https://sketchfab.com/models/b354421b86df4871b142e76d153026f2/embed?transparent=1}{\includegraphics[width=.95\linewidth]{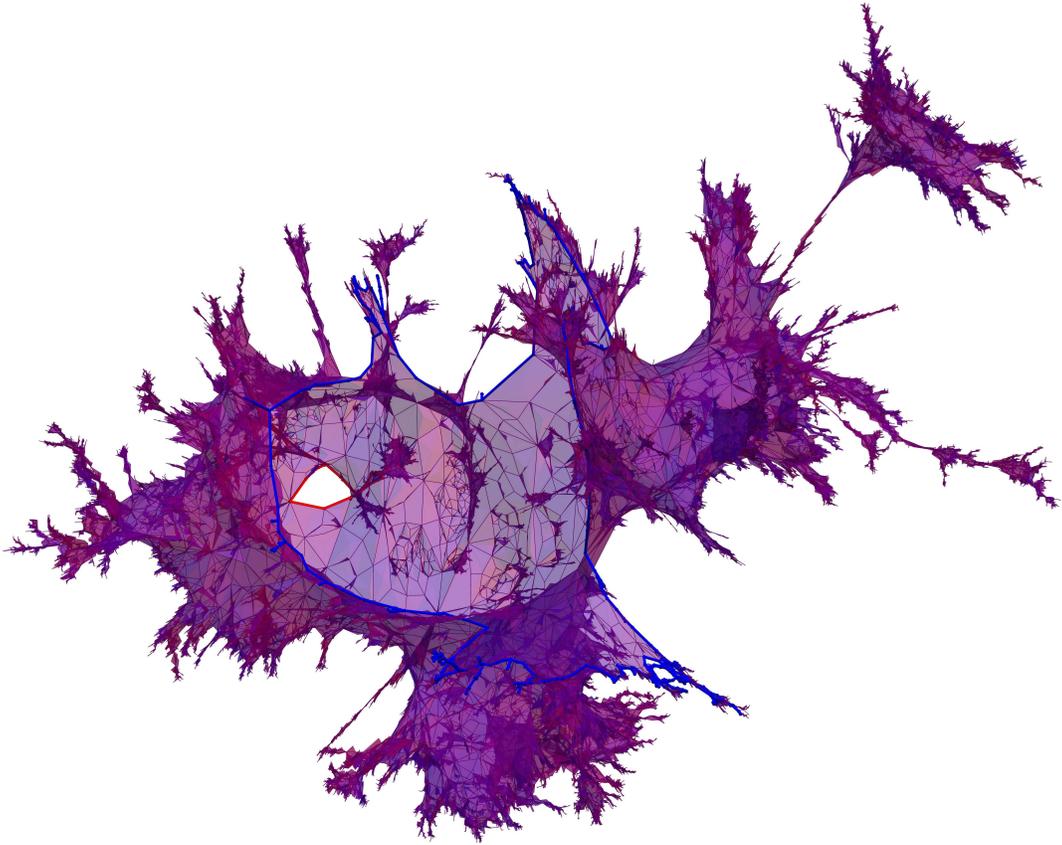}}
	\caption{Simulation of a uniformly sampled quadrangulation from $\Q_{50\,000,(200,300)}$ in genus~$2$. The $2$ boundary components are highlighted. \emph{Click on the picture or go to \nolinkurl{www.normalesup.org/}\texttildelow\nolinkurl{bettinel} in order to zoom or rotate the map.}}
	\label{simu}
\end{figure}

\section{Introduction}

The purpose of this work is twofold. We first construct a class of metric spaces we call \emph{Brownian surfaces}, which arise as the scaling limit of random discrete maps on general surfaces. They naturally generalize the famous \emph{Brownian map} \cite{legall11ubm,miermont11bms}, which corresponds to the case where the surface in question is the sphere. 

In a second step, we study all the geodesics starting from a ``uniformly'' chosen random point in a fixed Brownian surface. Our method allows to recover the main result of~\cite{legall08glp} and provides analogous results for general Brownian surfaces. As the topology of the spaces we consider here is richer than the topology of the sphere, we find new geometrical characterizations in terms of boundary points and pairs of geodesics whose concatenation forms a loop that is not homotopic to~$0$.

\begin{figure}[ht!]
	\centering\includegraphics{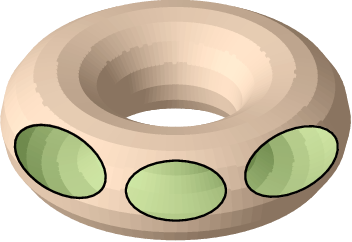}
	\caption{The surface with a boundary $\mathbf\Sigma_{1,3}^{\partial}$.}
	\label{s13}
\end{figure}

Recall that a \emph{surface with a boundary} is a non empty Hausdorff topological space in which every point has an open neighborhood homeomorphic to some open subset of~$\R\times\R_+$. Its \emph{boundary} is the $1$-dimensional manifold consisting of the points having a neighborhood homeomorphic to a neighborhood of $(0,0)$ in~$\R\times\R_+$. Note that, in particular, a surface without boundary is a surface with a boundary whose boundary is empty. In this work, we will only consider compact connected orientable surfaces with a boundary. By the classification theorem, they are characterized up to homeomorphisms by two nonnegative integers, the genus~$g$ and the number~$p$ of connected components of the boundary. We denote by~$\Sgp$ the unique (up to homeomorphism) compact orientable surface of genus~$g$ with~$p$ boundary components; it can be obtained from the compact orientable surface of genus~$g$ by removing~$p$ disjoint open disks whose boundaries are pairwise disjoint circles. See Figure~\ref{s13}.
\nomenclature[0]{$\Sgp$}{compact orientable surface of genus~$g$ with~$p$ boundary components}%

A map is a cellular embedding of a finite graph (possibly with multiple edges and loops) into a compact connected orientable surface without boundary,
considered up to orientation-preserving homeomorphisms. Cellular means that the faces of the map (the connected components of the complement of edges) are open $2$-cells, that is, homeomorphic to $2$-dimensional open disks. We will consider in this work \emph{maps with a boundary}, which means that some distinguished faces of the map are considered as ``holes.'' Given such a map, we may remove from every distinguished face an open disk whose boundary is a circle and obtain a surface with a boundary. It would probably be more satisfactory to directly define these objects as embedded into a surface with a boundary in such a way that the distinguished faces are homeomorphic to disks with an open disk removed and all other faces are $2$-cells but the previous definition happens to be more convenient to our purpose.

The natural problem of \emph{scaling limits} of random maps has generated many studies in the last decade. The most natural setting is the following. We consider maps as metric spaces, endowed with their natural graph metric. We choose uniformly at random a map of ``size''~$n$ in some class, rescale the metric by the proper factor, and look at the limit in the sense of the Gromov--Hausdorff topology~\cite{gromov99msr}. The size considered is often the number of faces of the map. From this point of view, the most studied class is the class of planar quadrangulations. The pioneering work of Chassaing and Schaeffer~\cite{chassaing04rpl} revealed that the proper scaling factor in this case is~$n^{-1/4}$. The problem was first addressed by Marckert and Mokkadem~\cite{marckert06limit}, who constructed a candidate limiting space called the \textit{Brownian map}, and showed the convergence toward it in another sense. Le~Gall~\cite{legall07tss} then showed the relative compactness of this sequence of metric spaces and that any of its accumulation points was almost surely of Hausdorff dimension~$4$. It is only recently that the solution of the problem was completed independently by Miermont~\cite{miermont11bms} and Le~Gall~\cite{legall11ubm}, who showed that the scaling limit is indeed the Brownian map. This last step, however, is not mandatory in order to identify the topology of the limit: Le~Gall and Paulin~\cite{legall08slb}, and later Miermont~\cite{miermont08sphericity}, showed that any possible limit is homeomorphic to the $2$-dimensional sphere. 

This line of reasoning lead the way to several extensions. The first kind of extension is to consider other classes of planar maps. Actually, Le~Gall already considered in~\cite{legall07tss} the classes of $\kappa$-angulations, for even $\kappa\ge 4$. In~\cite{legall11ubm}, he considered the classes of $\kappa$-angulations for $\kappa=3$ and for even $\kappa\ge 4$ as well as the case of Boltzmann distributions on bipartite planar maps, conditioned on their number of vertices. Another extension is due to Addario-Berry and Albenque~\cite{addarioalbenque2013simple}; in this work, they consider simple triangulations and simple quadrangulations, that is, triangulations and quadrangulations without loops and multiple edges. Together with Jacob and Miermont~\cite{BeJaMi14}, we later added the case of maps conditioned on their number of edges, and Abraham~\cite{abr13} considered the case of bipartite maps conditioned on their number of edges. In all these cases, the limiting space is always the same Brownian map (up to a multiplicative constant): we say that the Brownian map is \emph{universal} and we expect it to arise as the scaling limit of a lot more of natural classes of maps. A peculiar extension is due to Le~Gall and Miermont~\cite{legall09scaling} who consider maps with large faces, forcing the limit to fall out of this universality class: they obtain so-called \emph{stable maps}, which are related to stable processes.

Another kind of extension is to consider quadrangulations on a fixed surface that is no longer the sphere. The case of orientable surfaces of positive genus was the focus of~\cite{bettinelli10slr,bettinelli12tsl} and the case of the disk was considered in~\cite{bettinelli11slr}. In both these cases, the last step is still missing at the moment: the convergence is only known to hold along subsequences. In the present work, we complete this picture by considering quadrangulations on arbitrary surfaces with a boundary.

The starting point of these problems is a powerful bijective encoding of the maps in the studied class by simpler objects. In the case of planar quadrangulations, the bijection in question is the so-called Cori--Vauquelin--Schaeffer bijection~\cite{cori81planar,schaeffer98cac,chassaing04rpl} and the simpler objects are trees whose vertices carry integer labels satisfying some conditions. In the other cases, variants of this bijection are used \cite{bouttier04pml,chapuy07brm,poulscha06,ambjornbudd} and the encoding objects usually have a more intricate combinatorial structure.

At the present time, little is known about the limiting spaces. The principal result about the metric properties of the Brownian map is due to Le~Gall~\cite{legall08glp}. In this reference, Le~Gall studies all the geodesics from a uniformly chosen random point in the Brownian map. In particular, he shows that, toward a typical point, there is only one such geodesic and that for particular points there are up to three distinct such geodesics. Moreover, if we consider the unique geodesics toward two typical points, it is shown that they share a common part of positive length. We believe that Le~Gall's approach could be adapted to our framework of general Brownian surfaces. We also think that, in the particular case $p=0$, the results could be derived from Le~Gall's, by using a continuous analog to Chapuy's bijection~\cite{chapuy08sum} (see also \cite[Section~6]{bettinelli12tsl} for such an application). In the present work, however, we choose to use a third approach.

Let $\bp$, $x_1$, $x_2$, \ldots\ be fixed uniform points in a Brownian surface. Our strategy relies on the following facts, which hold almost surely.
\begin{enumerate}[($a$)]
	\item There is only one geodesic from~$\bp$ to~$x_i$ (Proposition~\ref{geodtyp}).\label{idea1}
	\item It is possible to ``entrap'' all the geodesics from~$\bp$ to any given point with geodesics from~$\bp$ to~$x_i$'s. (Section~\ref{secgengeod}).\label{idea2}
\end{enumerate}
Combining these facts, we are able to identify all the geodesics starting from~$\bp$. Fact~($\ref{idea1}$) will be shown by a technique inspired from~\cite{miermont09trm}. In this reference, a bijection between quadrangulations with~$q$ distinguished vertices and maps with~$q$ faces is provided. Very roughly speaking, a proper use of the parameters of this bijection in the particular case $q=2$ gives a parametrization of the geodesics between the two distinguished points. Fact ($\ref{idea2}$) is the core of our approach and will be explained in more details in the next section, as we need some more notation.

\bigskip

Throughout this work, we fix an integer $g\ge 0$ and work in fixed genus~$g$.
\nomenclature[11]{$g\ge0$}{fixed genus in which we work}%
In order to ease the reading, we provide page~\pageref{secnot} a recapitulation of our most used notation. We will see during Section~\ref{secsl} that a Brownian surface may be constructed both as a quotient of $[0,1]$ and as a quotient of a more complicated underlying structure. As often as possible, $r$, $s$, $t$ will denote numbers in $[0,1]$, $a$, $b$, $c$ will denote points in the underlying structure and $x$, $y$, $z$ will denote  points in the Brownian surface.

\paragraph*{Acknowledgments.}

The author warmly thanks Gr\'egory Miermont for interesting discussions leading to the realization of this work. The author also thanks Curtis McMullen for suggesting the terminology of Brownian surfaces, as well as an anonymous referee for tediously reading a previous version of this work and for making numerous improving comments.

\section{Main results}
\subsection{Brownian surfaces}

Let~$\m$ be a map. We let $V(\m)$ denote its set of vertices, $E(\m)$ its set of edges, and $\vec E(\m)$ its set of half-edges (oriented edges).
\nomenclature[1m1]{$\m$, $\q$}{map, quadrangulation}%
\nomenclature[1m2]{$V$, $E$, $\vec E$}{vertex set, edge set, half-edge set of a map}%
We say that a face~$f$ is \emph{incident} to a half-edge~$e$ (or that~$e$ is incident to~$f$) if~$e$ is included in the boundary of~$f$ and is oriented in such a way that~$f$ lies to its left. The number of half-edges incident to a face is called its \emph{degree}. For any half-edge~$e$, we denote by~$\bar e$ its reverse, as well as~$e^-$ and~$e^+$ its origin and end.
\nomenclature[1me]{$e^-$, $e^+$, $\bar e$}{origin, end, reverse of the half-edge~$e$}%
We denote by~$d_\m$ the graph metric on~$\m$ defined as follows: for any $v,v'\in V(\m)$, the distance $d_\m(v,v')$ is the number of edges of any shortest path linking~$v$ to~$v'$.
\nomenclature[1m5]{$d_\m$}{graph metric on~$\m$}%

\begin{figure}[ht!]
		\psfrag{i}[][][.8]{$h_1$}
		\psfrag{j}[][][.8]{$h_2$}
		\psfrag{k}[][][.8]{$h_3$}
	\centering\includegraphics[width=88mm]{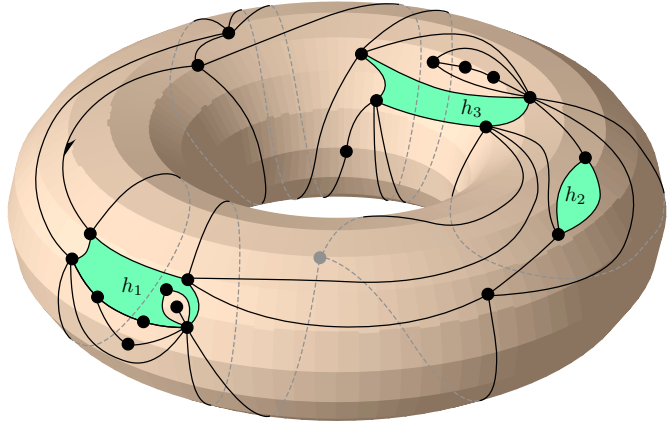}
	\caption{A quadrangulation from $\Q_{19,(4,1,2)}$ in genus~$1$. The half arrowhead symbolizes the root.}
	\label{exmap}
\end{figure}

For combinatorial reasons, we need to restrict ourselves to bipartite maps: a map is called \emph{bipartite} if its vertex set can be partitioned into two subsets such that no edge links two vertices of the same subset.  For $p\ge 0$, a \emph{quadrangulation with~$p$ boundary components} is a bipartite map having~$p$ distinguished faces~$h_1$, \dots, $h_p$ and whose other faces are all of degree~$4$. The distinguished faces will be called \emph{external faces} or \emph{holes}. The other faces will be called \emph{internal faces}. For technical reasons, the maps we consider will always implicitly be rooted, in the sense that they come with a distinguished half-edge, usually denoted by~$e_*$. The \emph{genus} of a map is defined as the genus of the surface into which it is embedded. For $n\in\Zp$ and $\si=(\si^1,\dots,\si^p) \in \N^p$ (with the convention that $\N^0\de \{\emptyset\}$), we define the set~$\Qns$ of all genus~$g$ quadrangulations with~$p$ boundary components having~$n$ internal faces and such that~$h_i$ is of degree $2\si^i$, for $1\le i \le p$.
\nomenclature[1q2]{$e_*$}{root half-edge}%
\nomenclature[3d1]{$p \ge 0$}{number of holes}%
\nomenclature[3d2]{$h_1$, \dots, $h_p$}{holes in~$\q$ or in~$\m$}%
\nomenclature[3d3]{$n \ge 0$}{number of internal faces of~$\q$ / number of edges not incident to a hole in~$\m$}%
\nomenclature[3d4]{$\si_i$}{half boundary length of $h_i$ in~$\q$ / boundary length of $h_i$ in~$\m$}%
\nomenclature[3d5]{$\Qns$}{set of genus~$g$ quadrangulations with~$p$ boundary components having~$n$ internal faces and such that~$h_i$ is of degree $2\si^i$, for $1\le i \le p$}%
See Figure~\ref{exmap} for an example and Figure~\ref{simu} for a computer simulation of a large random quadrangulation.

The Gromov--Hausdorff distance between two compact metric spaces $(\X,\delta)$ and $(\X',\delta')$ is defined by
$$\dGH\big((\X,\delta),(\X',\delta')\big) \de \inf \big\{ \dish\big(\varphi(\X),\varphi'(\X')\big)\!\big\},$$
where the infimum is taken over all isometric embeddings $\varphi : \X \to \X''$ and $\varphi':\X'\to \X''$ of~$\X$ and~$\X'$ into the same metric space $(\X'', \delta'')$, and $\dish$ stands for the usual Hausdorff distance between compact subsets of $\X''$. This defines a metric on the set of isometry classes of compact metric spaces (\cite[Theorem~7.3.30]{burago01cmg}), making it a Polish space (this is a simple consequence of \cite[Theorem~7.4.15]{burago01cmg}).

The following result unifies and generalizes the main results from \cite{legall07tss,legall08slb,miermont08sphericity,bettinelli10slr,bettinelli12tsl, bettinelli11slr}.

\begin{thm}\label{cvthm}
Let us consider an integer $p\ge 0$, positive real numbers $\si_\8^1$, \dots, $\si_\8^p>0$ and a sequence of $p$-uples $\si_n=(\si_n^1,\dots,\si_n^p) \in \N^p$ such that $\si_n^i/{\sqrt{2n}} \to \si^i_\8$, for $1\le i \le p$. Let $\q_n$ be uniformly distributed over~$\Qnsn$. Then, from any increasing sequence of integers, we may extract a subsequence $(n_k)_{k\ge 0}$ such that there exists a random metric space $(\qis,\disig)$ satisfying
$$\Big( V(\q_{n_k}),\lp\frac 9 {8n_k}\rp^{1/4} d_{\q_{n_k}} \Big) \tolk \big(\qis,\disig\big)$$
in the sense of the Gromov--Hausdorff topology.

Moreover, regardless of the choice of the sequence of integers, the limiting space $(\qis,\disig)$ is almost surely homeomorphic to~$\Sgp$, has Hausdorff dimension~$4$, and every of the~$p$ connected components of its boundary has Hausdorff dimension~$2$.
\end{thm}

\begin{rem}
The constant $(8/9)^{1/4}$ is irrelevant in this statement. We chose to let it figure for the sake of consistency with previous works and because of following definitions. This constant is inherent to the case of quadrangulations: we believe that the same statement should hold with the same limiting spaces for other classes of maps embedded in the same surface and satisfying mild conditions, up to modifying this constant.
\end{rem}

In the spherical case ($(g,p)=(0,0)$), these results were shown by Le~Gall \cite{legall07tss} and the topology was identified by Le~Gall and Paulin \cite{legall08slb} (and also later by Miermont~\cite{miermont08sphericity}). Using two different approaches, Le~Gall~\cite{legall11ubm} and Miermont~\cite{miermont11bms} later substantially improved this result by showing that the extraction is not necessary in this case. It is strongly believed that it is not necessary in the general case either; this is the main focus of the works in progress~\cite{BeMi15I,BeMi15II}.
Theorem~\ref{cvthm} in the particular case $p=0$, $g>0$ can be found in \cite{bettinelli10slr,bettinelli12tsl} and the case $p=1$, $g=0$ is treated in~\cite{bettinelli11slr}.

Notice that, although it seems reasonable that the limiting space will have genus at most~$g$ and at most~$p$ holes, it is not clear a priori that it will be homeomorphic to~$\Sgp$. We could imagine that some handles ``disappear'' or that some holes ``merge'' into a single hole. This does not happen; loosely speaking, this means that a uniform quadrangulation is sufficiently well spread over the surface, taking a macroscopic amount of space inbetween the holes and on every handle. Another noticeable fact is that the boundary of every hole is homeomorphic to a circle whereas, in the discrete picture, the holes do not in general have a simple curve as a boundary.

We believe that it is possible to extend this result to the case where some~$\sigma_\8^i$'s are equal to~$0$. More precisely, making the extra hypothesis that the extraction in Theorem~\ref{cvthm} is not necessary, we conjecture the following.
\begin{conj}
Let $0\le p'\le p$ be integers, let $\si_\8^1$, \dots, $\si_\8^{p'}>0$ be positive real numbers and set $\si_\8^{p'+1}=0$, \dots, $\si_\8^p=0$. Let also $\si_n=(\si_n^1,\dots,\si_n^p) \in \N^p$ be $p$-uples such that $\si_n^i/{\sqrt{2n}} \to \si^i_\8$, for $1\le i \le p$ and consider a random variable~$\q_n$ uniformly distributed over~$\Qnsn$. Then, the following convergence holds in the sense of the Gromov--Hausdorff topology:
$$\Big( V(\q_{n}),\lp\frac 9 {8n}\rp^{1/4} d_{\q_{n}} \Big) \tol \big(\qi^{\si'},\di^{\si'}\big),$$
where $\big(\qi^{\si'},\di^{\si'}\big)$ is the limiting space of Theorem~\ref{cvthm} corresponding to $\si_\8^1$, \dots, $\si_\8^{p'}$.
\end{conj}

We proved this conjecture in~\cite{bettinelli11slr} in the case $g=0$, $p'=0$ and $p=1$ (where it is known \cite{legall11ubm,miermont11bms} that the extraction is not necessary) and we think that with some substantial work, our method can be adapted to this framework, provided the uniqueness of the limiting space. We decided, however, not to pursue this in the present paper as we prefer to concentrate on geodesics in the limiting spaces.

\subsection{Geodesics in Brownian surfaces}\label{secgeodbs}

We now turn to the geodesics in a fixed Brownian surface. We fix a subsequence $(n_k)_{k\ge 0}$ along which the convergence of Theorem~\ref{cvthm} holds and we consider the corresponding Brownian surface $(\qis,\disig)$. Recall that, in a compact metric space $(\X,\delta)$, a \emph{geodesic} from a point~$x\in\X$ to a point~$y\in \X$ is a continuous path $\wp: [0,\delta(x,y)] \to \X$ such that $\wp(0)=x$, $\wp(\delta(x,y))=y$ and $\delta(\wp(s),\wp(t))=|t-s|$ for every $s,t\in[0,\delta(x,y)]$. The space $(\X,\delta)$ is called a \emph{geodesic space} if any two points are connected by at least one geodesic. Because the metric space associated with a map is at $\dGH$-distance at most~$1/2$ from a geodesic space and because the Gromov--Hausdorff limit of a sequence of geodesic spaces is also a geodesic space (see \cite[Theorem~7.5.1]{burago01cmg}), we see from Theorem~\ref{cvthm} that any Brownian surface is a geodesic space. Although we do not know how to characterize all the geodesics in a Brownian surface, we are able to describe all the geodesics starting from a distinguished point. Before stating our results, we still need to introduce some more notions.

\begin{figure}[ht!]
		\psfrag{i}[][]{$h_1$}
		\psfrag{j}[][]{$h_2$}
		\psfrag{s}[][]{$s$}
		\psfrag{t}[][]{$t$}
		\psfrag{u}[][]{$s'$}
		\psfrag{v}[][]{$t'$}
		\psfrag{r}[][]{\textcolor{gray}{$\bp$}}
		\psfrag{p}[][]{\textcolor{red!90!black}{$\Phi_s$}}
		\psfrag{q}[][]{\textcolor{red!90!black}{$\Phi_t$}}
		\psfrag{m}[][]{\textcolor{red!90!black}{$\Phi_{s'}$}}
		\psfrag{n}[][]{\textcolor{red!90!black}{$\Phi_{t'}$}}
	\centering\includegraphics[width=88mm]{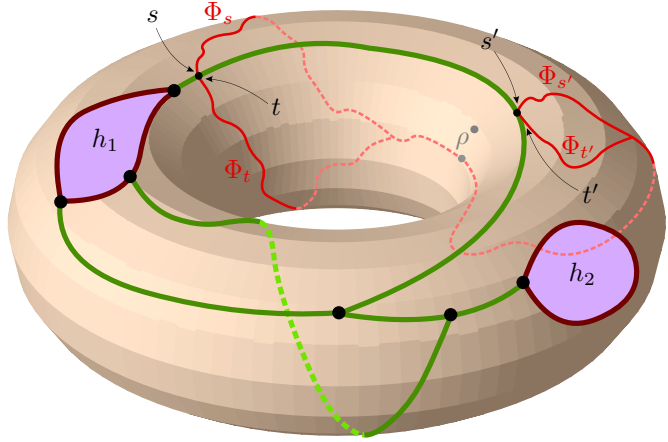}
	\caption{Description of a Brownian torus with~$2$ holes. The backbone is the thick green and burgundy structure. It is shaped like a map with~$2$ holes and one face but we should think of the edges as having some ``length.'' Imagine that, on both sides of the edges, independent Brownian forests are embedded, except inside the~$2$ holes. In other words, Brownian CRT's are grafted according to some Poissonian distribution. The union of these forests and the backbone constitutes the underlying structure~$\Mi$.}
	\label{scheme}
\end{figure}

We will see later that the Brownian surface $(\qis,\disig)$ may be constructed as a quotient of an underlying structure~$\Mi$. This structure consists of a \emph{backbone}, which roughly captures its homotopy type, on which Brownian forests are glued. It generalizes Aldous's CRT, which is the analog in the spherical case. The rigorous definition of this object will be given in Section~\ref{secquot}, after some grounds will have been established. At this stage, we only give an informal description and refer to Figure~\ref{scheme} for a visual support. For the moment, do not pay attention to $s$, $t$, $s'$, $t'$, and the corresponding paths $\Phi$ on the figure; these elements will be of interest in a short while. 

It is also important to know at this point that~$\Mi$ comes with a natural distinguished point noted~$\bp$. This is the base point for the geodesics we will study. In general, there may exist simple loops in $(\qis,\disig)$ that are not homotopic to~$0$. This motivates the following definition. We let $\NC(\bp,\qis)$ be the set of points~$x\in\qis$ for which there exists at least a pair $\{\wp,\wp'\}$ of geodesics from~$\bp$ to~$x$ such that the concatenation of~$\wp$ with the reverse of~$\wp'$ is not homotopic to~$0$. The Hausdorff dimension of this set can be computed:

\begin{thm}\label{thmdimnc}
The set $\NC(\bp,\qis)$ is a.s.\ empty if $g=0$ and $p\in\{0,1\}$; it a.s.\ has Hausdorff dimension~$2$ otherwise.
\end{thm}

We define the \emph{order} of a point $a\in\Mi$ as the maximal number of connected components of $\Mi\bs\{a\}$ restricted to the neighborhoods of~$a$. Basic properties of Brownian forests show that this number is a most~$3$. The backbone $\cB\subseteq\Mi$ is the set of points~$a$ either of order~$2$ and such that no connected components of $\Mi\bs\{a\}$ are real trees or of order~$3$ and such that at most one connected component of $\Mi\bs\{a\}$ is a real tree. For instance, note that $\cB=\emptyset$ in the case $(g,p)=(0,0)$ of the sphere, as~$\Mi$ is a real tree in this case; in fact, this is the only case in which $\cB=\emptyset$. Proposition~\ref{propid} shows in particular that only points of order~$1$ can be identified in the quotient, so that we may define the order of a point in~$\qis$ as the order of any of its representatives in~$\Mi$. Note that this definition a priori depends on~$\Mi$. We denote by~$\partial\qis$ the boundary of~$(\qis,\disig)$. The following theorem generalizes the main result of~\cite{legall08glp}:

\begin{thm}\label{thmgeod}
The following properties hold almost surely.
\begin{enumerate}[($i$)]
	\item For all $x\in\qis$, the number of distinct geodesics from~$\bp$ to~$x$ is equal to the order of~$x$ minus $\un{\{x\in\partial \qis\}}$. In particular, this number is typically~$1$, at most~$2$ for the boundary points and at most~$3$ for the interior points.\label{thmgeodi}
	\item The canonical projection of~$\cB$ in~$\qis$ is the union of~$\partial\qis$ and the set $\NC(\bp,\qis)$.\label{thmgeodii}
	\item The set of points reachable by~$3$ distinct geodesics and for which every pair $\{\wp,\wp'\}$ of geodesics is such that the concatenation of~$\wp$ with the reverse of~$\wp'$ is not homotopic to~$0$ is finite: its cardinality~$H$ is 
	\begin{itemize}
		\item equal to~$0$ if $g=0$ and $p\in\{0,1\}$;
		\item equal to $4g-2$ if $g\ge 1$ and $p=0$;
		\item a random variable whose distribution only depends on~$g$, $p$ and~$\sigma$ (it can even be computed) and whose support is $\{0,1,\dots,4g+p-2\}$ otherwise.\label{thmgeodiii}
	\end{itemize}
	\item We suppose that $p\neq 0$ in order to avoid trivialities. The set $\partial\qis \cap \NC(\bp,\qis)$ is finite: its cardinality is equal to $4g+2p-2-H$.\label{thmgeodiv}
\end{enumerate}
\end{thm}

Let us first informally describe in more details what happens on the example of Figure~\ref{scheme}. The boundary~$\partial\qis$ corresponds to the burgundy edges and the set $\NC(\bp,\qis)$ corresponds to the green edges. The projection in~$\qis$ of the vertices of the 3-face map corresponding to the backbone constitutes a finite set of cardinality $4g+2p-2$ (by a simple application of the Euler characteristic formula). Among these points, $H$ are not incident to any hole and $4g+2p-2-H$ are incident to some hole ($H=4$ on the picture). As the $3$-face map of the picture is actually a random variable, we see that~$H$ is also random and the exact distribution can be computed once we know the distribution of this $3$-face map. Moreover, these $4g+2p-2$ points are precisely the ones described in~($\ref{thmgeodiii}$) or~($\ref{thmgeodiv}$), depending on whether they belong to~$\partial\qis$ or not.

We see that, in particular, our notion of order for a point in~$\qis$ only depends on~$\bp$. It is possible to define a uniform measure on~$\qis$, as the projection of the Lebesgue measure on~$[0,1]$ (it can also be seen as the limit in some sense of the uniform measure on $V(\q_n)$). Typical in~($\ref{thmgeodi}$) means that, if a point is distributed according to this measure, the number for this point is a.s.~$1$. We can also show that~$\bp$ is uniformly distributed over~$\qis$ and the ``invariance under uniform rerooting'' property (see \cite[Section~8]{legall08glp} for precise statements in the spherical case, which can easily be extended to general Brownian surfaces) yields a similar result if we replace~$\bp$ with another uniformly chosen random point. Another important fact to notice is that the distribution of~$H$ does not depend on the subsequence $(n_k)_{k\ge 0}$. This reinforces the conjecture that the extraction in Theorem~\ref{cvthm} is actually not necessary.

Using basic properties of~$\Mi$, we can classify the points of~$\qis$ as follows. The set of points reachable from~$\bp$ by a unique geodesic a.s.\ has Hausdorff dimension~$4$. The set of points reachable from~$\bp$ by at least two distinct geodesics is a dense subset of~$\qis$, which is relatively small (we can show that its Hausdorff dimension is a.s.~$2$). An important part of these points moreover belong to $\NC(\bp,\qis)$, which is also a set of dimension a.s.~$2$. Finally, the set of points that can be reached by three geodesics is countable and only a finite number of these points are such that the concatenation of any two geodesics is not homotopic to~$0$.


Note also that, as a random map on a surface is locally planar (at least far from the boundary), it is to be expected that, in the vicinity of~$\bp$, every Brownian surface should look alike. In this regard, it is probable that~($\ref{thmgeodi}$) in the latter vicinity can also be directly derived from Le~Gall's results~\cite{legall08glp} by an absolute continuity argument.

\medskip

Let us comment a little more about the idea behind the proof of Theorem~\ref{thmgeod}. It comes from a complete description of the geodesics from~$\bp$: we show that, in the language of~\cite{legall08glp}, all the geodesics from~$\bp$ are \emph{simple} (see Proposition~\ref{propallgeod}). We can picture~$\Mi$ as a continuous analog to a discrete map with~$1$ internal face and~$p$ holes. As for a real tree (which is the analog to a discrete tree, that is, a map with~$1$ face in genus~$0$), we will see that it is possible to define a contour order around the ``internal face'' of~$\Mi$. The order of a point~$a\in\Mi$ is then the number of times~$a$ is visited in this contour order plus~$1$ if~$a$ is incident to a hole (because the holes are not visited). To each time~$a$ is visited in the contour order correspond a simple geodesic. Furthermore, Lemma~\ref{lemhomot} gives a criterion to decide whether the concatenation of two simple geodesics aiming at the same point are homotopic to~$0$ or not, in terms of the visiting times of the corresponding point in~$\Mi$.

Going back to Figure~\ref{scheme}, we denote by~$\Phi_s$ the simple geodesic corresponding to the time~$s$. While performing the contour order, at time~$s$, we are visiting some forest (even possibly its floor) grafted on a half-edge of the $3$-face map: we denote this half-edge by $\tilde {e_s}$. We consider a point in~$\Mi$ that is visited both at times~$s$ and~$t$. Then the criterion of Lemma~\ref{lemhomot} states that the concatenation of~$\Phi_s$ with the reverse of~$\Phi_t$ is homotopic to~$0$ if and only if $\tilde {e_{s}} = \tilde {e_{t}}$. For instance, on the picture of Figure~\ref{scheme}, $\tilde {e_{s}} \neq \tilde {e_{t}}$ so that~$\Phi_s$ and~$\Phi_t$ constitute a loop that is not homotopic to~$0$. On the contrary, $\Phi_{s'}$ and~$\Phi_{t'}$ constitute a loop that is homotopic to~$0$. This explains the correspondence between the sets considered in Theorem~\ref{thmgeod} and the elements of the $3$-face map. Finally, the $-\un{\{x\in\partial \qis\}}$ in~($\ref{thmgeodi}$) simply comes from the fact that there are no geodesics inside the holes.

We can now give the intuition of Fact~($\ref{idea2}$) from the introduction. For simplicity, let us suppose that~$a$ is a point of order~$1$ and let~$s$ be the time at which it is visited in the contour order. By density, we can find times~$s_i$ and $s_j$ as close as wanted to~$s$ such that $s_i < s <s_j$ and such that the corresponding points~$x_i$ and~$x_j$ in~$\qis$ are from our fixed sequence of uniform points. By Fact~($\ref{idea1}$), the geodesics from~$x_i$ and~$x_j$ to~$\bp$ are necessarily simple geodesics. Using the description of simple geodesics, we can show that these geodesics have to merge within the vicinity of the point~$x$ corresponding to~$a$ in~$\qis$. As a result, every geodesic from~$x$ to~$\bp$ is ``trapped'' between these two geodesics and can only escape through the segment of~$\Mi$ connecting the points corresponding to~$s_i$ and~$s_j$. Using fine properties of the labels on~$\Mi$, we show that it is possible to select~$s_i$ and~$s_j$ in such a way that the geodesic does not escape through this segment. Hence it has to coincide with the other two geodesics, at least after they merge.

\bigskip

To complete this presentation, we also recover the same ``confluence property'' as the one obtained by Le~Gall in the spherical case \cite[Corollary~7.7]{legall08glp}.
\begin{prop}[{Corollary~\ref{corconfgeod}}]
Almost surely, for every $\eps>0$, there exists $\eta\in\,(0,\eps)$ such that all the geodesics from~$\bp$ to points outside of the ball of radius~$\eps$ centered at~$\bp$ share a common initial part of length~$\eta$.
\end{prop}

To end this section, let us mention that, as observed by Le~Gall in the spherical case, our results share surprising similarities with the results of Myers~\cite{myers36} in the context of differential geometry, where the surfaces considered posses by far more regularities. In the latter reference, it is shown that the locus of points reachable from a fixed point by at least two geodesics in a surface without boundary is a linear graph. Moreover, the subtraction of this locus from the surface yields a unique open $2$-cell and the number of distinct geodesics to a point of the locus is equal to the order of this point in the locus (where the definition of order matches our definition of order in~$\Mi$).

\subsection{Geodesics in large quadrangulations}\label{secgeodlq}

Theorem~\ref{thmgeod} allows us to derive asymptotic results for large quadrangulations. We adapt the presentation and technique of Le~Gall \cite{legall08glp} to obtain them. Note that, in this section, we do not restrict our attention to subsequences along which the convergence of Theorem~\ref{cvthm} holds.

We define the distance between a path~$\wp$ going successively through the vertices $v_1$, \ldots, $v_{k}$ and a path~$\wp'$  going successively through the vertices $v'_1$, \ldots, $v'_{k'}$ in the same map~$\m$ by
$$\dpath(\wp,\wp')=\max_{i\ge 0}\big\{d_{\m}(v_{i\wedge k},v'_{i\wedge k'})\big\}.$$
For $\eps>0$ and $v$, $v'\in V(\q_n)$, we let $\Geo^\eps_n(v,v')$ be the set of paths from~$v$ to~$v'$ of length at most $(1+\eps)\,d_{q_n}(v,v')$. For $\delta>0$, we also let $\Mult_{n,\delta}^\eps(v,v')$ be the maximal integer~$r$ such that there exist~$r$ paths $\wp_1$, \ldots, $\wp_r\in \Geo^\eps_n(v,v')$ satisfying $\dpath(\wp_i,\wp_j)\ge \delta n^{1/4}$ whenever $i\neq j$.

Let $(\eps_n)_n$ be a fixed deterministic sequence of nonnegative real numbers such that $\eps_n\to 0$ as $n\to\8$ and, conditionally given~$\q_n$, let $v_n^\ooo$ be a uniformly distributed vertex in~$V(\q_n)$. We denote by $\partial V(\q_n)$ the set of vertices in $V(\q_n)$ that are incident to a hole. We start with the translation of~($\ref{thmgeodi}$) in terms of large discrete maps.

\begin{prop}\label{propdiscgeodi}
Let $v_n$ denote a uniformly distributed vertex in~$V(\q_n)$, independent from~$v_n^\ooo$. Then, for every $\delta>0$, $\Pb\big(\Mult_{n,\delta}^{\eps_n}(v_n^\ooo,v_n)=1\big) \to 1$ as $n\to \8$.

Moreover, for every $\delta>0$ and for every sequence $(\eta_n)_n$ of positive numbers such that $\eta_n\to 0$ as $n\to\8$, $\Pb\big(\exists v\in V(\q_n) : d_{\q_n}(v,\partial V(\q_n))\le n^{1/4}\eta_n,\,\Mult_{n,\delta}^{\eps_n}(v_n^\ooo,v)\ge 3\big) \to 0$  as $n\to \8$.

Finally, for every $\delta>0$, $\Pb\big(\exists v\in V(\q_n) : \Mult_{n,\delta}^{\eps_n}(v_n^\ooo,v)\ge 4\big) \to 0$  as $n\to \8$ and
$$\lim_{\delta\to 0} \liminf_{n\to\8}  \Pb\big(\exists v\in V(\q_n): \Mult_{n,\delta}^{0}(v_n^\ooo,v)= 3\big) = 1.$$
\end{prop}

Note that we chose to state the result with geodesics in the displayed equation, instead of ``approximate'' geodesics because it is actually stronger. Of course, the weaker result with~$\eps_n$ in place of~$0$ also holds. The first part of Proposition~\ref{propdiscgeodi} roughly means that, in a large uniform quadrangulation, there is essentially a unique macroscopic geodesic from~$v_n^\ooo$ to a uniformly chosen vertex and this remains true even if we allow some ``slack'' to the geodesics. The second part states that the points asymptotically close to the boundary cannot be reached from~$v_n^\ooo$ by more than two distinct macroscopic geodesics, even if ``slack'' is allowed. The last part says that asymptotically, there exist points that can be reached from~$v_n^\ooo$ by three distinct macroscopic geodesics but no points can be reached by four distinct macroscopic geodesics, even if ``slack'' is allowed.

\bigskip

We now complete the interpretation. For $\eps>0$ and $v$, $v'\in V(\q_n)$, we let $\HMult_{n}^\eps(v,v')$ be the maximal integer~$r$ such that there exist~$r$ paths $\wp_1$, \ldots, $\wp_r\in \Geo^\eps_n(v,v')$ such that, whenever $i\neq j$, the concatenation of~$\wp_i$ with the reverse of~$\wp_j$ is not homotopic to~$0$. Here, homotopic to~$0$ means homotopic to~$0$ in the surface with a boundary corresponding to the map, that is, the surface in which the map is embedded, with open disks removed from the external faces. Note that, as the topology is preserved at the limit (handles and holes do not vanish), we can find a (random) $\delta$ small enough such that, for large~$n$, $\HMult_{n}^\eps(v,v')\le \Mult_{n,\delta}^\eps(v,v')$. With a little bit of work, this yields the following statement:

\begin{prop}\label{prophmult}
We have $\Pb\big(\exists v\in V(\q_n) : \HMult_{n}^{\eps_n}(v_n^\ooo,v)\ge 4\big) \to 0$ as $n\to\8$.
\end{prop}

In order to give a discrete version of~($\ref{thmgeodii}$), we need to know that the encoding bijection gives rise to a particular subset $\cB(\q_n,v_n^\ooo)\subseteq V(\q_n)$, which is the discrete counterpart of~$\cB$. For the reader already familiar with this kind of bijections, it is the $2$-core of the encoding object, that is, the set of vertices remaining after iteratively removing from it all the vertices of degree~$1$. See the end of Section~\ref{sec1pt} for a precise definition. 

It is not hard to see that all the points $v\in\cB(\q_n,v_n^\ooo)\bs\partial V(\q_n)$ satisfy $\HMult_{n}^{0}(v_n^\ooo,v)\ge 2$. The following proposition roughly states that all such points should be close to this set.

\begin{prop}\label{propdiscgeodii}
For every $\delta>0$,
\begin{align*}
\Pb\big(\exists v\in V(\q_n)\, :\, d_{\q_n}(v,\cB(\q_n,v_n^\ooo)\bs\partial V(\q_n))\ge \delta n^{1/4},\, \HMult_{n}^{\eps_n}(v_n^\ooo,v)\ge 2\big) \ton 0.
\end{align*}
\end{prop}

The following statement means that the vertices $v\in V(\q_n)$ such that $\HMult_{n}^{\eps_n}(v_n^\ooo,v)=3$ are essentially dispatched in small separate ``zones'' and that the number of these zones is asymptotically the random variable~$H$.

\begin{prop}\label{propdiscgeodiii}
Let $A_{n,\delta}^\eps(j)$ denote the event that there exist~$j$ vertices $v_1$, \ldots, $v_j\in V(\q_n)$ such that $\HMult_{n}^{\eps}(v_n^\ooo,v_i)=3$ for $1\le i \le j$ and $d_{\q_n}(v_i,v_{i'})\ge \delta n^{1/4}$ whenever $i\neq i'$. Then
$$\lim_{\delta\to 0}\liminf_{n\to\8}\Pb\big(A_{n,\delta}^{\eps_n}(j)\bs A_{n,\delta}^{\eps_n}(j+1) \big)= 
	\lim_{\delta\to 0}\limsup_{n\to\8}\Pb\big(A_{n,\delta}^{\eps_n}(j)\bs A_{n,\delta}^{\eps_n}(j+1) \big)= \Pb(H=j),$$
where~$H$ is the random variable defined in~($\ref{thmgeodiii}$).
\end{prop}

Finally, we have a similar result for the boundary vertices~$v$ for which $\HMult_{n}^{\eps_n}(v_n^\ooo,v)\ge 2$.

\begin{prop}\label{propdiscgeodiv}
We suppose that $p\neq 0$ and we denote by $B_{n,\delta}^\eps(j)$ the event that there exist~$j$ vertices $v_1$, \ldots, $v_j\in  \partial V(\q_n)$ such that $\HMult_{n}^{\eps}(v_n^\ooo,v_i)\ge 2$ for $1\le i \le j$ and $d_{\q_n}(v_i,v_{i'})\ge \delta n^{1/4}$ whenever $i\neq i'$. Then
$$\lim_{\delta\to 0}\liminf_{n\to\8}\Pb\big(B_{n,\delta}^{\eps_n}(j)\bs B_{n,\delta}^{\eps_n}(j+1) \big)=\Pb(H=4g+2p-2-j)$$
and the same statement holds if we replace the $\liminf$ with a $\limsup$.
\end{prop}

\subsection{Organization of the paper}

In Section~\ref{secenc}, we explain how to encode quadrangulations by simpler maps. We begin with the mapping for pointed quadrangulations, which constitutes the starting point of our study and we then present the two-point mapping, which will be crucially used in Section~\ref{secuniq}. Section~\ref{secsl} is devoted to scaling limits and the proof of Theorem~\ref{cvthm}. As it uses the general framework of \cite{bettinelli10slr,bettinelli12tsl,bettinelli11slr}, we only sketch it and point to the particularity of the cases considered here. We then study typical geodesics and show Fact~($\ref{idea1}$) in Section~\ref{sectypgeod} and finally, we study the general geodesics in Section~\ref{secgengeod} by the entrapment technique (Fact~($\ref{idea2}$)). This leads to Proposition~\ref{propallgeod}, which identifies all the geodesics from~$\bp$. Section~\ref{secremproofs} is devoted to the remaining proofs.

\section{Encoding quadrangulations}\label{secenc}

Let us now present the bijections allowing to encode quadrangulations carrying one or two distinguished vertices. Our description is a slight reformulation of the \BDG bijection~\cite{bouttier04pml} and the Miermont bijection~\cite{miermont09trm}. We refer the reader to these references for proofs. Throughout this section, $n$ and~$p$ denote nonnegative integers and $\si=(\si^1,\dots,\si^p)$ a $p$-uple of positive integers.

\subsection{Quadrangulations carrying one distinguished vertex}\label{sec1pt}

\begin{figure}[ht!]
		\psfrag{d}[][][.8]{$d$}
		\psfrag{e}[r][r][.8]{$d+1$}
		\psfrag{f}[l][l][.8]{$d+2$}
		\psfrag{g}[l][l][.8]{$d+1$}
		\psfrag{1}[][][.8]{$3$}
		\psfrag{2}[][][.8]{$4$}
		\psfrag{3}[][][.8]{$5$}
		\psfrag{4}[][][.8]{$6$}
		\psfrag{5}[][][.8]{$7$}
		\psfrag{a}[][][.8]{$c_i^0$}
		\psfrag{b}[l][l][.8]{$c_i^1$}
		\psfrag{c}[][][.8]{$c_i^2$}
		\psfrag{h}[][][.8]{$c_i^3$}
		\psfrag{i}[][][.8]{$c_i^4$}
		\psfrag{j}[][][.8]{$c_i^5$}
		\psfrag{k}[][][.8]{$c_i^6$}
	\centering\includegraphics{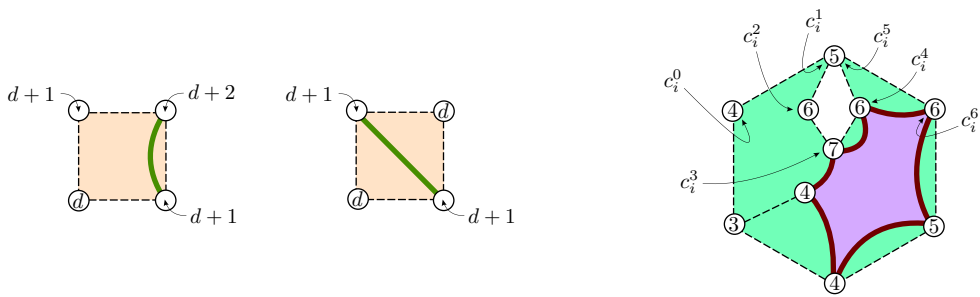}
	\caption{\textbf{Left.} Adding a new edge inside an internal face. \textbf{Right.} Adding~$\si^i$ new edges inside the external face~$h_i$. On this example, $\si^i=6$.}
	\label{newedges}
\end{figure}

Let $\q\in \Qns$ be a quadrangulation and $v^\ooo \in V(\q)$ one of its vertices. We assign labels to the vertices of~$\q$ as follows: for every vertex $v\in V(\q)$, we set $\lab(v)\de d_\q(v^\bullet,v)$. Because~$\q$ is by definition bipartite, the labels of both ends of any edge differ by exactly~$1$. As a result, the internal faces can be of two types: the labels around the face are either $d$, $d+1$, $d+2$, $d+1$, or $d$, $d+1$, $d$, $d+1$ for some~$d$. We add a new edge inside every internal face following the convention depicted on the left part of Figure~\ref{newedges}.

A \emph{corner} is an angular sector delimited by two successive half-edges in the contour of a face. The vertex located at the end of the first half-edge is called the vertex \emph{incident} to the corner. If~$c$ is a corner incident to a vertex~$v$, we write $\lab(c)\de\lab(v)$ with a slight abuse of notation. For each~$i$, we let $c_i^0$, $c_i^1$, \dots, $c_i^{2\si^i-1}$ be the corners of~$h_i$ read in clockwise order, starting at an arbitrary corner (and we use the convention $c_i^{2\si^i}\de c_i^0$). We link together in a cycle the corners~$c_i^k$'s such that $\lab(c_i^{k+1})=\lab(c_i^k)-1$, as shown on the right part of Figure~\ref{newedges}. Note that, because $\lab(c_i^{k+1})-\lab(c_i^k) \in \{-1,+1\}$, there are exactly~$\si^i$ such corners.

We then only keep the new edges we added and the vertices in $V(\q)\bs\{v^\ooo\}$. The object we obtain is a labeled map~$\m$ of genus~$g$ with $p+1$ faces. There is an obvious correspondence between the external faces of~$\q$ and~$p$ of the faces of~$\m$. Let~$h_1$, \dots, $h_p$ also denote these faces in~$\m$. Note that, by construction, these faces all have a simple boundary and are of degrees~$\si^1$, \dots, $\si^p$. Remark also that~$v^\ooo$ lies within the remaining face of~$\m$, which we denote by~$f^\ooo$.

\begin{figure}[ht!]
		\psfrag{f}[][][.8]{$f^\ooo$}
		\psfrag{v}[][][.8]{$v^\ooo$}
		\psfrag{q}[][]{$\q$}
		\psfrag{i}[][][.8]{$h_1$}
		\psfrag{j}[][][.8]{$h_2$}
		\psfrag{k}[][][.8]{$h_3$}
		\psfrag{0}[][][.7]{$0$}
		\psfrag{1}[][][.7]{$1$}
		\psfrag{2}[][][.7]{$2$}
		\psfrag{3}[][][.7]{$3$}
		\psfrag{4}[][][.7]{$4$}
	\centering\includegraphics[width=.95\linewidth]{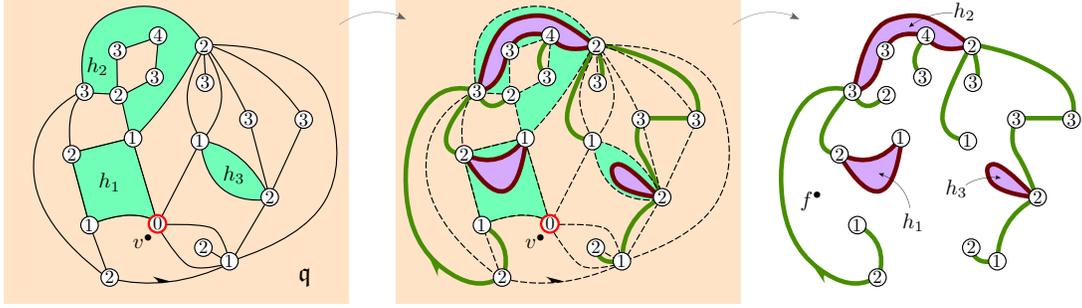}
	\caption{The construction for a map in $\Q_{12,(2,4,1)}$ in the case $g=0$.}
	\label{bijqm}
\end{figure}

\begin{figure}[ht!]
		\psfrag{f}[][][.8]{$f$}
		\psfrag{d}[][][.8]{$d$}
		\psfrag{g}[l][l][.8]{$d+1$}
	\centering\includegraphics{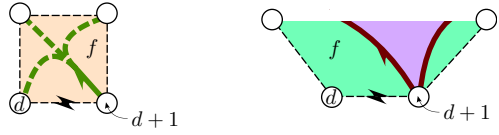}
	\caption{Rooting the map~$\m$. On the picture, the two possible roots for~$\q$ yielding the same root for~$\m$ are shown.}
	\label{rooting}
\end{figure}

We root the map~$\m$ as follows. Let~$e$ be the only half-edge among the root of~$\q$ and its reverse such that $\lab(e^+)=\lab(e^-)+1$, and let~$f$ be the face of~$\q$ that is incident to~$e$. If~$f$ is an internal face, the root of~$\m$ is the half-edge corresponding to the edge we added in~$f$, directed from~$e^+$. If~$f$ is an external face, there are two new half-edges inside~$f$ starting from~$e^+$; the root of~$\m$ is the one incident to~$f^\ooo$. See Figure~\ref{rooting}.

For each~$i$, let $\vh_i^{1}$, $\vh_i^2$, \dots, $\vh_i^{\si^i}$ be the half-edges incident to~$h_i$ in~$\m$, read in counterclockwise order around it. The labels of~$\m$ satisfy the following:
\begin{itemize}
	\item for $1\le i \le p$ and $1 \le j \le \si^i$, we have $\lab(\vh_i^{j -}) - \lab(\vh_i^{j +}) \ge -1$;	
	\item for any half-edge~$e\in \vec E(\m)$ such that neither~$e$ nor its reverse~$\bar e$ is incident to an~$h_i$, we have $|\lab(e^+) - \lab(e^-)|\le 1$.
\end{itemize}

We will consider the labels of~$\m$ up to an additive constant: we write
$$[\lab]\de \{v\in V(\m) \mapsto \lab(v)+a:\ a\in \Z\}$$
\nomenclature[2b5]{$(\m,[\lab])$}{labeled map}%
the class of~$\lab$ for this equivalence relation. We say that two faces are \emph{adjacent} if there exist a half-edge incident to one and whose reverse is incident to the other. Let $\Mnsun$ denote the set of genus~$g$ maps having $n+|\si|$ edges (where we write $|\si|\de\sum_{i=1}^p \si^i$) and $p+1$ faces denoted by $h_1$, \dots, $h_p$, $f^\ooo$ such that, for all~$i$, $h_i$ has a simple boundary, is of degree~$\si^i$ and is not adjacent to any other~$h_j$, and such that the root is not incident to any hole~$h_j$.
\nomenclature[4d1]{$\Mnsun$}{set of genus~$g$ maps having $n+$\textbar$\si$\textbar{} edges and $p+1$ faces denoted by $h_1$, \dots, $h_p$, $f^\ooo$ such that, for all~$i$, $h_i$ has a simple boundary, is of degree~$\si^i$ and is not adjacent to any other~$h_j$, and such that the root is not incident to any hole~$h_j$}%
Note that any edge is forbidden to be incident to two different holes, but there may exist vertices that are incident to two or more holes. We also denote by~$\LMnsun$ the set of pairs whose first coordinate lies in~$\Mnsun$ and whose second coordinate is an equivalence class of labeling functions on the vertices of the map satisfying the two previous itemized conditions.
\nomenclature[4d2]{$\LMnsun$}{set of pairs $(\m,[\lab])$ where $\m\in\Mnsun$ and~$\lab$ is a suitable labeling function}%

\begin{prop}\label{bijqvml}
The mapping $(\q,v^\ooo) \mapsto (\m,[\lab])$ is a two-to-one mapping from the set of quadrangulations in~$\Qns$ carrying one distinguished vertex to the set $\LMnsun$.
\end{prop}

At this point, we may properly define the set $\cB(\q,v^\ooo)$ appearing in Section~\ref{secgeodlq}. Let $(\m,[\lab])$ denote the map corresponding to~$(\q,v^\ooo)$ by the previous mapping. Then $\cB(\q,v^\ooo)$ is the $2$-core of~$\m$, that is, the set of vertices of~$\m$ remaining after iteratively removing from~$\m$ all its vertices of degree~$1$.

\subsection{Quadrangulations carrying two distinguished vertices}\label{sectwopoint}

In the case of doubly-pointed quadrangulations, the construction is very similar. We consider a quadrangulation $\q\in \Qns$ and $v^\ooo$, $v^{\ooo\ooo} \in V(\q)$ two of its vertices such that $d\de d_{\q}(v^\ooo,v^{\ooo\ooo})\ge 2$.
\nomenclature[2b1]{$v^\ooo$, $v^{\ooo\ooo}$}{distinguished vertices in~$\q$}%
We also need an integer $\lambda\in \{1,2,\dots, d-1\}$.
We assign labels to the vertices of~$\q$ in a similar fashion as in the last section but we now take into account the extra vertex~$v^{\ooo\ooo}$: we define
\begin{equation}\label{eqdefl}
\lab(v)\de \min \big( d_\q(v^\ooo,v), d_\q(v^{\ooo\ooo},v) + 2\lambda- d \big) ,\qquad v\in V(\q).
\end{equation}
A simple way to picture this function is to imagine water flowing through the edges of~$\q$ at rate~$1$, starting from two sources: one located at the vertex~$v^\ooo$ opened at time~$0$ and one located at the vertex~$v^{\ooo\ooo}$ opened at time $2\lambda-d$.

We then follow rigorously the same procedure as in last section, and we obtain a class of labeled maps $(\m,[\lab])$ lying in $\LMnsde$, which is defined exactly as $\LMnsun$, with the exception that the map now has $p+2$ faces $h_1$, \dots, $h_p$, $f^\ooo$, $f^{\ooo\ooo}$ instead of~$p+1$.
\nomenclature[2b3]{$f^{\ooo}$, $f^{\ooo\ooo}$}{face of~$\m$ carrying~$v^\ooo$, face of~$\m$ carrying~$v^{\ooo\ooo}$}%
\nomenclature[4d4]{$\LMnsde$}{similar definition as~$\LMnsun$ with one additional face~$f^{\ooo\ooo}$}%
It turns out that~$v^\ooo$ and~$v^{\ooo\ooo}$ do not belong to the same face of~$\m$ and we let~$f^\ooo$ be the face containing~$v^\ooo$ and~$f^{\ooo\ooo}$ be the one containing~$v^{\ooo\ooo}$. We also define the analog~$\Mnsde$ of~$\Mnsun$ with $p+2$ faces.
\nomenclature[4d3]{$\Mnsde$}{similar definition as~$\Mnsun$ with one additional face~$f^{\ooo\ooo}$}%
Let us emphasize at this point that, because the holes have a simple boundary and are not adjacent with one another, every edge of~$\m$ is incident to at least one of the faces~$f^\ooo$ or~$f^{\ooo\ooo}$.

\begin{figure}[ht!]
		\psfrag{v}[][][.8]{$v^\ooo$}
		\psfrag{w}[][][.8]{$v^{\ooo\ooo}$}
		\psfrag{q}[][]{$\q$}
		\psfrag{i}[][][.8]{$h_1$}
		\psfrag{j}[][][.8]{$h_2$}
		\psfrag{k}[][][.8]{$h_3$}
		\psfrag{f}[][][.8]{$f^\ooo$}
		\psfrag{g}[][][.8]{$f^{\ooo\ooo}$}
		\psfrag{0}[][][.7]{$0$}
		\psfrag{1}[][][.7]{$1$}
		\psfrag{2}[][][.7]{$2$}
		\psfrag{3}[][][.7]{$3$}
		\psfrag{4}[][][.7]{$4$}
		\psfrag{9}[][][.7]{-$1$}
		\psfrag{8}[][][.7]{-$2$}
	\centering\includegraphics[width=.95\linewidth]{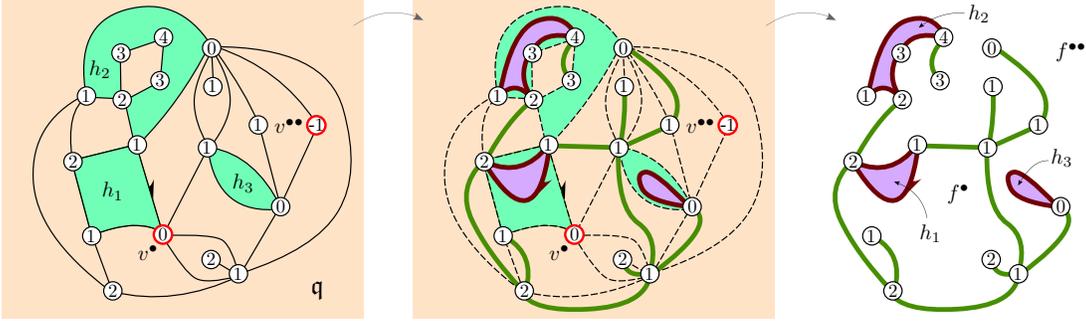}
	\caption{On this example, $d=3$ and $\lambda=1$. It can be checked for example that the vertices corresponding in~$\q$ to the three vertices labeled~$1$ on the interface between~$f^\ooo$ and~$f^{\ooo\ooo}$ are at distance~$1$ from~$v^\ooo$ and at distance~$2$ from~$v^{\ooo\ooo}$.}
	\label{bijqm2}
\end{figure}

\begin{prop}
This construction provides a two-to-one mapping from
$$\Big\{ \big(\q,(v^\ooo, v^{\ooo\ooo}),\lambda\big):\, \q\in \Qns,\, (v^\ooo, v^{\ooo\ooo}) \in V(\q)^2\ \big|\ d_{\q}(v^\ooo,v^{\ooo\ooo})\ge 2,\, 1\le\lambda\le d_{\q}(v^\ooo,v^{\ooo\ooo})-1 \Big\}$$
to $\LMnsde$.
\end{prop}

\begin{rem}
In~\cite{miermont09trm}, Miermont considers what he calls \emph{delay vectors}. In the particular case of $2$-pointed quadrangulations, it is more natural to use the integer~$\lambda$ instead. The construction we give here corresponds to the delay vector $[0,2\lambda-d]$ with his notation. It can moreover be observed that, if~$u$ is any of the vertices where the water from the two sources meet for the first time, then $\lambda=d_\q(v^\ooo,u)$ and $d-\lambda=d_\q(v^{\ooo\ooo},u)$. In some sense, the integer~$\lambda$ parametrizes the geodesics between~$v^\ooo$ and~$v^{\ooo\ooo}$. This observation will be the key point of the proof of Proposition~\ref{geodtyp}. 
\end{rem}

\subsection{Reverse mappings}\label{secrevmap}

In both cases, the reverse mapping uses the same procedure. We describe it briefly in this section. Let $(\m,[\lab])$ be in either $\LMnsun$ or $\LMnsde$. First, we add inside~$f^\ooo$ a new vertex~$v^\ooo$ with label
$$\lab(v^\ooo) \de  \min_{u\in f^{\ooo}} \lab(u)-1,$$
where the notation $u\in f$ means that~$u$ belongs to the boundary of~$f$. Following the counterclockwise order around~$f^\ooo$, we draw arcs linking every corner to the first subsequent corner that has a strictly smaller label. If no such corners exist, which means that the corner we are visiting has a label that is minimal on~$f^\ooo$, we draw the arc from the corner to the extra vertex~$v^\ooo$. It is possible to draw these arcs in such a way that they do not cross each other or the edges of~$\m$. In the case where $(\m,[\lab])\in \LMnsde$, we follow the same procedure for~$f^{\ooo\ooo}$: we add a vertex~$v^{\ooo\ooo}$ inside it with label $\lab(v^{\ooo\ooo}) \de \min_{u\in f^{\ooo\ooo}} \lab(u)-1$ and we draw arcs in the same way.
 
When removing the edges of~$\m$, we are left with a quadrangulation~$\q$. Each hole~$h_i$ naturally corresponds to an external face of~$\q$, which we also denote by~$h_i$. The root of~$\q$ is defined as the arc drawn from the corner preceding the root of~$\m$, oriented in one direction or the other.

If $(\m,[\lab])\in \LMnsun$, its two pre-images are the pairs $(\q,v^\ooo)$, where~$\q$ is rooted in the two possible ways. If $(\m,[\lab])\in \LMnsde$, its two pre-images are
$$\big( \q,(v^\ooo, v^{\ooo\ooo}),\lambda\big)\qquad\text{ where }\lambda\de\frac12 \big( \lab(v^{\ooo\ooo})-\lab(v^{\ooo}) + d_\q(v^\ooo, v^{\ooo\ooo}) \big),$$
with the two possible rootings of~$\q$.

\begin{lem}\label{lemlambda}
In the case of the two-point mapping, the following holds:
\begin{enumerate}[($i$)]
	\item if~$v\in f^\ooo$, then $\lab(v)-\lab(v^\ooo)=d_\q(v^\ooo,v)$ and $\lab(v)-\lab(v^{\ooo\ooo})\le d_\q(v^{\ooo\ooo},v)$;
	\item if~$v\in f^{\ooo\ooo}$, then $\lab(v)-\lab(v^\ooo)\le d_\q(v^\ooo,v)$ and $\lab(v)-\lab(v^{\ooo\ooo})= d_\q(v^{\ooo\ooo},v)$;
	\item if~$v\in f^{\ooo\ooo}$, then $d_\q(v^\ooo,v) \ge \lambda$, if $v\in f^{\ooo}$, then $d_\q(v^{\ooo\ooo},v) \ge d_\q(v^{\ooo\ooo},v^\ooo)-\lambda$, and
	$$\lambda = \min_{v\in f^\ooo\cap f^{\ooo\ooo}} d_\q(v^\ooo,v)=\min_{v\in f^\ooo\cap f^{\ooo\ooo}} \lab(v)-\min_{v\in f^\ooo} \lab(v) +1.$$
\end{enumerate}
\end{lem}

\begin{pre}
($i$) and ($ii$). By~\eqref{eqdefl}, for any vertex~$v$, we have $\lab(v)-\lab(v^\ooo)\le d_\q(v^\ooo,v)$ and $\lab(v)-\lab(v^{\ooo\ooo})\le d_\q(v^{\ooo\ooo},v)$ (note that the representative of~$[\lab]$ appearing in~\eqref{eqdefl} is the one vanishing at~$v^\ooo$). Moreover, if $v\in f^\ooo$, the successive arcs drawn from~$v$ to~$v^\ooo$ form a path in~$\q$ of length $\lab(v)-\lab(v^\ooo)$, so that $d_\q(v^\ooo,v)\le\lab(v)-\lab(v^\ooo)$. This gives ($i$) and ($ii$) is proven in a similar fashion.

($iii$). Let us fix $v\in f^{\ooo\ooo}$. By the triangle inequality and ($ii$), 
$$d_\q(v^\ooo, v^{\ooo\ooo})\le d_\q(v^\ooo,v) + \lab(v)-\lab(v^{\ooo\ooo}) \le 2d_\q(v^\ooo,v) + \lab(v^\ooo) -\lab(v^{\ooo\ooo}).$$
The first part of the statement follows and the second part is obtained by the same argument. Finally, if~$v\in f^\ooo\cap f^{\ooo\ooo}$ is on a geodesic linking~$v^\ooo$ to~$v^{\ooo\ooo}$, all the previous inequalities are equalities and $d_\q(v^\ooo,v) = \lambda$. The first equality follows and the second is obtained thanks to ($i$). Note also that the condition of non adjacency of the holes makes it clear that $f^\ooo\cap f^{\ooo\ooo}\neq\emptyset$.
\end{pre}

\subsection{Further decomposition}

A labeled map $(\m,[\lab])$ from $\LMnsun$ or $\LMnsde$ can be decomposed into simpler objects: a \emph{scheme}, which in some sense accounts for the homotopy type of~$\m$, and a collection of forests indexed by some half-edges of the scheme. The labeling function naturally gives rise to labels on the vertices of these forests as well as to bridges recording the labels on the cycles of the map~$\m$.

\begin{rems}
The case $(\m,[\lab])\in \mM_{n,\emptyset}^1$ in genus~$0$ is somehow degenerate. Indeed, in this case, $\m$ is merely a plane tree and cannot be further decomposed in our sense. Until further notice, we suppose that we are not in this case.
\end{rems}
\begin{rems}
Up to a slight difference caused by the root, a scheme is sometimes called \emph{kernel} in graph theory. We chose to stick with the terminology of scheme, which seems more common in the context of maps.
\end{rems}

\subsubsection{Decomposition of the map}

Let us first focus on the map~$\m$ and keep the labels for later. We refer to Figure~\ref{decm} for visual support. We iteratively remove from~$\m$ all its vertices of degree~$1$ that are not extremities of its root~$e_*$. The set of vertices remaining at this point is called the \emph{floor} of~$\m$. Among these vertices, some are called \emph{nodes}: all vertices of degree~$3$ or more are nodes and,
\begin{itemize}
	\item if~$e_*^-$ is of degree~$1$, then~$e_*^-$ is a node;
	\item if~$e_*^+$ is of degree~$1$, then~$e_*^+$ is a node;
	\item if neither~$e_*^-$ nor~$e_*^+$ has degree~$1$, then~$e_*^-$ is a node.
\end{itemize}
On this map, the vertices that are not nodes are of degree~$2$ and are arranged into chains joining nodes. We define the map~$\s$ by replacing each of these chains by a single edge. The root of~$\s$ is defined as the edge replacing the chain that contains~$e_*$, oriented in the same direction as~$e_*$. The map~$\s$ is a scheme in the following sense.

\begin{defi}
A \emph{scheme} of type $(p,1)$ is a genus~$g$ map with $p+1$ faces denoted by $h_1$, \dots, $h_p$, $f^\ooo$, whose root is not incident to any~$h_j$, and that satisfies the following conditions. For all~$i$, $h_i$ has a simple boundary and is not adjacent to any~$h_j$. There may only be one vertex with degree~$1$ or~$2$: if it has degree~$1$, then it is an extremity of the root; if it has degree~$2$, then it is the origin of the root.

A scheme of type $(p,2)$ is a genus~$g$ map with $p+2$ faces denoted by $h_1$, \dots, $h_p$, $f^\ooo$, $f^{\ooo\ooo}$ and satisfying the same conditions.
\end{defi}
\nomenclature[5s1]{$\s$}{scheme}%

\begin{rem}
A more conventional definition would be to forbid vertices of degree less than~$2$. Our choice of ``keeping the root present'' in the scheme is done for combinatorial reasons. It is only in Section~\ref{secremproofs} that we will need to get rid of this ``root,'' which does not play any geometrical part. Compare the definition of the backbone $\cB(\q,v^\ooo)$ used in the geometrical statements of Section~\ref{secgeodlq} with the definition of the floor used for combinatorial purposes.
\end{rem}

\begin{figure}[ht!]
		\psfrag{i}[][][.8]{$h_1$}
		\psfrag{f}[][][.8]{$f^\ooo$}		
	\centering\includegraphics[width=.95\linewidth]{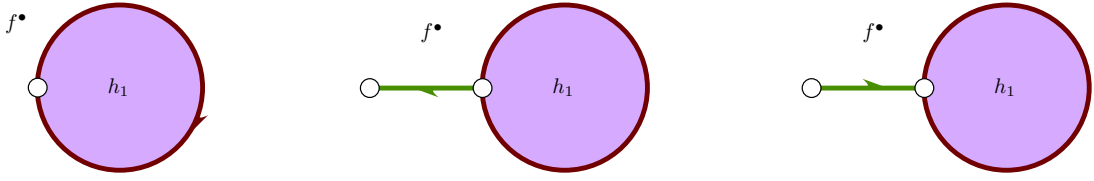}
	\caption{The three elements of $\Sg_{1,1}$ in genus $g=0$.}
	\label{sg11}
\end{figure}

Let~$\Sg_{p,1}$ and~$\Sg_{p,2}$ be the finite sets of schemes of type $(p,1)$ and $(p,2)$.
\nomenclature[5s3]{$\Sg_{p,q}$}{set of schemes of type $(p,q)$}%
For example, in genus $g=0$, the set $\Sg_{1,1}$ contains the three maps represented on Figure~\ref{sg11}.

We will use the following formalism for forests.
\begin{defi}
A \emph{forest} of length $\xi\ge 1$ and mass~$m\ge 0$ is an ordered family of~$\xi$ trees (which can be defined as planar one-faced maps) with total number of edges equal to~$m$. Let~$\F_\xi^m$ denote the set of these forests.
\end{defi}

It will be convenient to systematically add a $\xi+1$-th tree to a forest of~$\F_\xi^m$ consisting of one single vertex. The \emph{floor} of a forest is the set of the root vertices of its trees, with the extra vertex-tree included. In the drawings, we will also add extra edges linking the elements of the floor (as on Figure~\ref{decm}).

\begin{figure}[ht!]
		\psfrag{i}[][][.8]{$h_1$}
		\psfrag{j}[][][.8]{$h_2$}
		\psfrag{f}[][][.8]{$f^\ooo$}
		\psfrag{g}[][][.8]{$f^{\ooo\ooo}$}
		\psfrag{1}[][][.8]{$\vec H_1$}
		\psfrag{2}[][][.8]{$\vec H_2$}
		\psfrag{I}[][][.8]{$\vec I$}
		\psfrag{B}[][][.8]{$\vec B$}
	\centering\includegraphics[width=8cm]{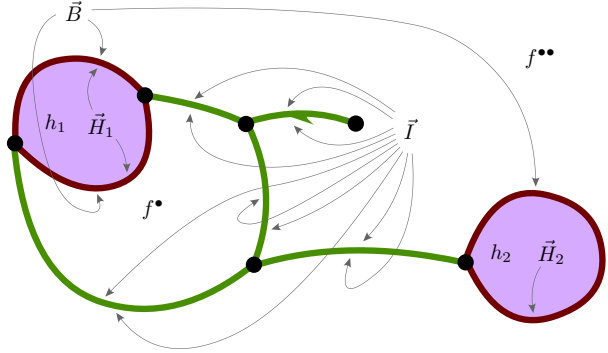}
	\caption{Notation for the different types of half-edges of a scheme. The set~$\vec F(\s)$ is the union of~$\vec B(\s)$ and~$\vec I(\s)$.}
	\label{notation}
\end{figure}

The half-edges of the scheme~$\s$ are of several different types. Let~$\vec H_i(\s)$ be the set of half-edges incident to the hole~$h_i$ and let~$\vec F(\s)$ be the set of half-edges that are not incident to any holes. We break the set $\vec F(\s)$ into two subsets: let $\vec I(\s)\de \{e\in \vec F(\s) : \bar e\in \vec F(\s)\}$ and $\vec B(\s)\de \{e\in \vec F(\s) : \bar e\notin \vec F(\s)\}$. See Figure~\ref{notation}. The letter $I$ stands for ``internal half-edges'' and $B$ stands for ``boundary half-edges.''
\nomenclature[5s11]{$\vec I(\s)$, $\vec B(\s)$}{internal half-edges of~$\s$, boundary half-edges of~$\s$}%

Every half-edge $e\in \vec F(\s)$ naturally corresponds to a forest~$\f^e$ defined as follows. Let $f\in\{f^\ooo,f^{\ooo\ooo}\}$ be the face incident to~$e$ and let $e'\in \vec F(\s)$ be the half-edge preceding~$e$ in the contour order of~$f$. By definition, the half-edges~$e$ and~$e'$ correspond to chains of half-edges in~$\m$: let~$\hat e$ and~$\hat e'$ be the last half-edges of these chains. The forest~$\f^e$ corresponds to the set of half-edges in~$\m$ visited between~$\hat e'$ and~$\hat e$ ($\hat e'$ excluded, $\hat e$ included) in the contour order of~$f$. See Figure~\ref{decm}.

For $e\in \vec F(\s)$, we denote by $\xi^e\ge 1$ and $m^e\ge 0$ the length and mass of the forest~$\f^e$: $\f^e\in\F_{\xi^e}^{m^e}$.
\nomenclature[5s5a]{$\xi^e\ge 1$}{length of the forest~$\f^e$}%
\nomenclature[5s5b]{$m^e\ge 0$}{mass of the forest~$\f^e$}%

\begin{prop}
\nomenclature[5s2]{$q\in\{1,2\}$}{number of distinguished vertices in the bijection}%
For $q\in\{1,2\}$, the above decomposition provides a bijection between the set~$\M_{n,\si}^q$ and the set of all pairs 
$$\big(\s,(\f^e)_{e\in\vec F(\s)}\big)$$
where $\s\in\Sg_{p,q}$ and such that there exist a collection of positive integers $(\xi^e)_{e\in\vec E(\s)}$ and a collection of nonnegative integers $(m^e)_{e\in\vec F(\s)}$ satisfying the following:
\begin{itemize}
	\item for all $e \in \vec F(\s)$, $\f^e\in\F_{\xi^e}^{m^e}$;
	\item for all $e \in \vec E(\s)$, $\xi^{\bar e}=\xi^e$;
	\item for $0\le i \le p$, $\sum_{e\in\vec H_i(\s)} \xi^{e}=\si^i$;
	\item $\sum_{e\in\vec F(\s)} m^e + \frac12 \sum_{e\in\vec E(\s)} \xi^{e} =n + |\si|$. 
\end{itemize}
\end{prop}

\subsubsection{Decomposition of the labeled map}

Let us now turn to the labels. The terminology we use here should become clear in a moment.

\begin{defi}
We call \emph{interior bridge} of length $\xi\ge 1$ from $a\in\Z$ to $b\in \Z$ a sequence of integers $(\bb(0), \dots, \bb(\xi))$ such that $\bb(0)=a$, $\bb(\xi)=b$ and, for all $0 \le i \le \xi -1$, we have $\bb(i+1) - \bb(i) \in \{-1,0,1\}$. We write $\fI_\xi(a,b)$ the set of these interior bridges.

We call \emph{boundary bridge} of length $\xi\ge 1$ from $a\in\Z$ to $b\in \Z$ a sequence $(\bb(0), \dots, \bb(\xi))$ of integers  such that $\bb(0)=a$, $\bb(\xi)=b$ and, for all $0 \le i \le \xi -1$, we have $\bb(i+1) - \bb(i) \ge -1$. We write $\fB_\xi(a,b)$ the set of these boundary bridges.
\end{defi}

\begin{defi}
A \emph{labeled forest} is a pair $(\f,\ell)$ where~$\f$ is a forest and $\ell:V(\f) \to \Z$ is a function satisfying the following:
\begin{itemize}
	\item for all~$u$ lying in the floor of~$\f$, $\ell(u)=0$;
	\item if~$u$ and~$v$ are linked by an edge, then $|\ell(u) - \ell(v)|\le 1$.
\end{itemize}
We denote by $\fF_\xi^m$ the set of labeled forests of length~$\xi$ and mass~$m$.
\end{defi}

\begin{figure}[ht!]
		\psfrag{i}[][][.8]{$h_1$}
		\psfrag{j}[][][.8]{$h_2$}
		\psfrag{f}[][][.8]{$f^\ooo$}
		\psfrag{g}[][][.8]{$f^{\ooo\ooo}$}
		\psfrag{m}[][][.8]{$(\m,[\lab])$}
		\psfrag{s}[][][.8]{$\s$}
		\psfrag{b}[][][.8]{$\bb^e$}
		\psfrag{l}[][][.8]{$(\f^e,\lab^e)$}
		\psfrag{0}[][][.7]{$0$}
		\psfrag{1}[][][.7]{$1$}
		\psfrag{2}[][][.7]{$2$}
		\psfrag{3}[][][.7]{$3$}
		\psfrag{4}[][][.7]{$4$}
		\psfrag{9}[][][.7]{-$1$}
		\psfrag{8}[][][.7]{-$2$}
		\psfrag{7}[][][.7]{-$3$}		
	\centering\includegraphics[width=.95\linewidth]{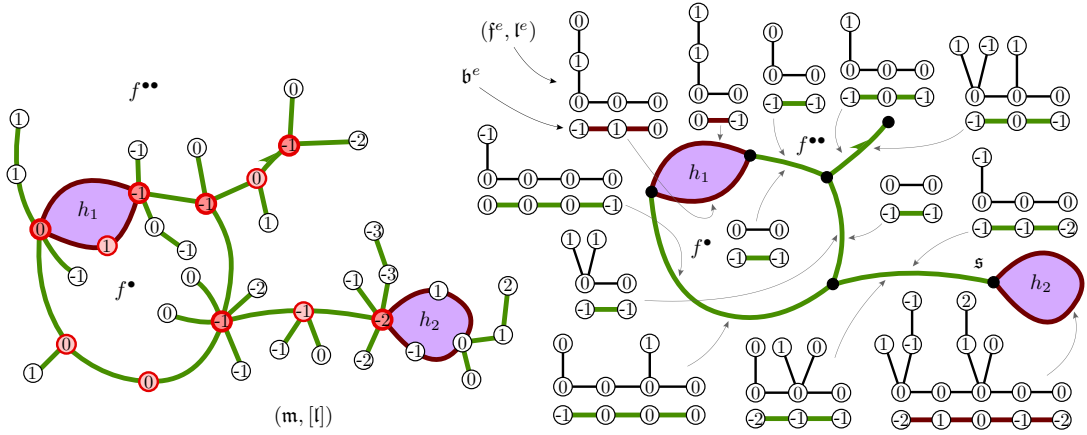}
	\caption{Decomposition of a labeled map from $\mM_{33,(3,4)}^2$ in genus~$0$ into a scheme $\s\in\Sg_{2,2}$, a collection of labeled forests $(\f^e,\lab^e)_{e\in\vec F(\s)}$ and a collection of bridges $(\bb^e)_{e\in\vec F(\s)}$. The bridges have different colors, depending on whether they are interior bridges or boundary bridges. On the left, the floor and the~$6$ nodes of~$\m$ are represented in red and with thicker outlines.}
	\label{decm}
\end{figure}

There is a trivial one-to-one correspondence between the nodes of~$\m$ and the vertices of~$\s$ so that~$\lab$ naturally gives a canonical labeling of the vertices of~$\s$ as follows. Let $v_*\in V(\s)$ be the origin of the root in~$\s$ and let $v\in V(\s)$. We denote by~$v'_*$ and $v'\in V(\m)$ the corresponding nodes in~$\m$ and we set $l^v \de \lab(v') - \lab(v_*')$. Let $e\in \vec F(\s)$. It naturally corresponds to a chain~$e_1$, $e_2$, \dots, $e_{\xi^e}$ of half-edges in~$\m$. We define the bridge
$$\bb^e \de \big(\lab(e_1^-) - \lab(v_*'),\lab(e_2^-) - \lab(v_*'),\dots,\lab(e_{\xi^e}^-) - \lab(v_*'),\lab(e_{\xi^e}^+) - \lab(v_*')\big).$$
\nomenclature[5s6]{$\bb^e$}{bridge corresponding to the half-edge~$e\in\s$}%
The constraints on~$\lab$ show that, if $e\in \vec I(\s)$ then $\bb^e \in \fI_{\xi^e}({l^{e^-}\! ,l^{e^+}})$ and, if $e\in \vec B(\s)$ then $\bb^e \in \fB_{\xi^e}({l^{e^-}\! ,l^{e^+}})$. Moreover, the forest~$\f^e$ from last section naturally inherits from~$\lab$ a labeling function~$\lab^e$ defined as follows. Let $u\in V(\f^e)$ and let $\rho\in V(\f^e)$ be the root of the tree to which~$u$ belongs. These vertices correspond to two vertices $u'$ and~$\rho' \in V(\m)$. We set $\lab^e(u)\de \lab(u')-\lab(\rho')$. See Figure~\ref{decm}.

\begin{prop}\label{propdec}
For $q\in\{1,2\}$, the above decomposition provides a bijection between the set~$\mM_{n,\si}^q$ and the set of all triples 
$$\big(\s,(\f^e,\lab^e)_{e\in\vec F(\s)},(\bb^e)_{e\in\vec F(\s)}\big)$$
\nomenclature[5s4]{$(\f^e,\lab^e)$}{labeled forest corresponding to the half-edge $e\in\s$}%
where $\s\in\Sg_{p,q}$ and such that there exist a collection of positive integers $(\xi^e)_{e\in\vec E(\s)}$, a collection of nonnegative integers $(m^e)_{e\in\vec F(\s)}$ and a collection of integers $(l^v)_{v\in V(\s)}$ satisfying the following:
\begin{itemize}
	\item for all $e \in \vec F(\s)$, $(\f^e,\lab^e)\in\fF_{\xi^e}^{m^e}$;
	\item for all $e \in \vec E(\s)$, $\xi^{\bar e}=\xi^e$;
	\item $l^{v_*}=0$, where~$v_*$ is the origin of the root of~$\s$;
	\item for all $e \in \vec I(\s)$, $\bb^e \in \fI_{\xi^e}({l^{e^-}\! ,l^{e^+}})$ and $\bb^{\bar e} = (\bb^e(\xi^e),\dots, \bb^e(0))$;
	\item for all $e \in \vec B(\s)$, $\bb^e \in \fB_{\xi^e}({l^{e^-}\! ,l^{e^+}})$;
	\item for $0\le i \le p$, $\sum_{e\in\vec H_i(\s)} \xi^{e}=\si^i$;
	\item $\sum_{e\in\vec F(\s)} m^e + \frac12 \sum_{e\in\vec E(\s)} \xi^{e} =n + |\si|$. 
\end{itemize}
\end{prop}

Note that the three collections of integers are entirely determined by the triple of scheme, forests and bridges.

\subsubsection{Encoding by real-valued functions}\label{secrealvf}

\begin{figure}[ht!]
	\begin{minipage}{0.55\linewidth}
		\psfrag{1}[][][.8]{$1$}
		\psfrag{2}[][][.8]{$2$}
		\psfrag{3}[][][.8]{$3$}
		\psfrag{4}[][][.8]{$4$}
		\psfrag{0}[][][.8]{$0$}
		\psfrag{9}[][][.8]{-$1$}
		\psfrag{8}[][][.8]{-$2$}
		\psfrag{a}[r][r][.8]{$\f^e(0)$, $\f^e(8)$}
		\psfrag{b}[r][r][.8]{$\f^e(1)$, $\f^e(5)$, $\f^e(7)$}
		\psfrag{c}[r][r][.8]{$\f^e(2)$, $\f^e(4)$}
		\psfrag{d}[r][r][.8]{$\f^e(3)$}
		\psfrag{e}[][][.8]{$\f^e(6)$}
		\psfrag{f}[][r][.8]{$\f^e(9)$}
		\psfrag{g}[][r][.8]{$\f^e(10)$}
	\centering\includegraphics[height=25mm]{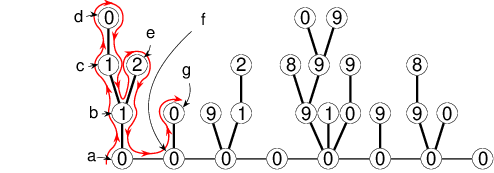}
	\end{minipage}
	\begin{minipage}{0.43\linewidth}
		\psfrag{C}[][][.8]{\textcolor{red}{$C^e$}}
		\psfrag{L}[][][.8]{\textcolor{blue}{$L^e$}}
	\centering\includegraphics[height=25mm]{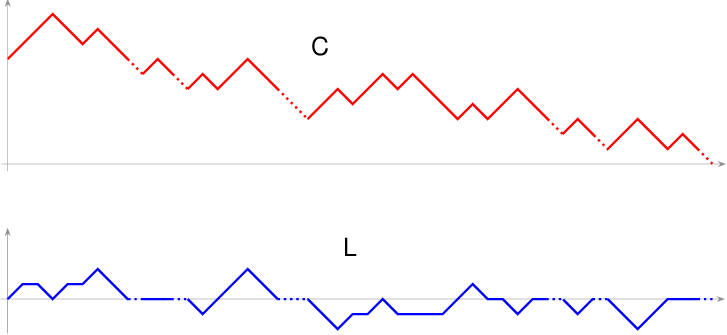}
	\end{minipage}
	\caption{The contour pair of a labeled forest from $\fF_7^{20}$. On the right, the paths are dashed on the intervals corresponding to edges linking elements of the floor.}
	\label{forest}
\end{figure}

For $e\in\vec F(\s)$, we encode the labeled forest $(\f^e,\lab^e)$ by its so-called \emph{contour pair} $(C^e,L^e)$ defined as follows. We see~$\f^e$ as a planar one-faced map with $2m^e+\xi^e$ edges (recall that we add a vertex-tree at the end and join by edges the elements of the floor). First, let $\f^e(0)$, $\f^e(1)$, \dots, $\f^e(2m^e+\xi^e)$ be the vertices of~$\f^e$ read in counterclockwise order around the face, starting at the first corner of the first tree. The \emph{contour function} $C^e:[0,2m^e+\xi^e]\to \R_+$ and the \emph{label function} $L^e:[0,2m^e+\xi^e]\to \R$ are defined by
$$C^e(i) \de d_{\f^e}\big(\f^e(i),\f^e(2m^e+\xi^e)\big) \sand L^e(i) \de \lab^e(\f^e(i)),\qquad 0 \le i \le 2m^e+\xi^e,$$
and linearly interpolated between integer values (see Figure~\ref{forest}).

We also define the function $B^e:[0,\xi^e]\to \R$ by
$$B^e(i) \de \bb^e(i),\qquad 0 \le i \le \xi^e,$$
and we linearly interpolate it between integer values. We will use the standard notation
$$\underline{X}(s) \de \inf_{0\le t\le s} X(t)$$
for the past infimum of any process~$X$. Remark that the function
\begin{equation*}\label{lele}
s \in [0, 2m^e+\xi^e] \mapsto L^e(s) + B^e\big( \xi^e - \underline C^e(s) \big)
\end{equation*}
records the labels up to an additive constant of the part of~$\m$ that corresponds to~$\f^e$. This will become useful in Section~\ref{secquot}.

\section{Scaling limit}\label{secsl}

From now on, we fix an integer $p\ge 0$ and a sequence $\si_n=(\si_n^1,\dots,\si_n^p) \in \N^p$ of $p$-uples such that
$$\si_{(n)}^i \de \frac{\si_n^i}{\sqrt{2n}} \to \si^i_\8 \in (0,\8),\qquad 1\le i \le p.$$
We set $\si_\8\de(\si_\8^1,\dots,\si_\8^p)$ and $\si_{(n)}\de \big(\si_{(n)}^1,\dots,\si_{(n)}^p\big)$ the rescaled version of~$\si_n$. Recall also that the genus~$g\ge 0$ is fixed. We furthermore suppose in this section that $(g,p)\neq (0,0)$ as this case requires special care.

\begin{rem}
Throughout this paper, the notation with a parenthesized~$n$ will always refer to suitably rescaled objects and the notation with an~$\8$ symbol will refer to limiting objects, as in the definitions above.
\end{rem}

Let~$\q_n$ be a random variable uniformly distributed over the set~$\Qnsn$ and let~$v_n^\ooo\in V(\q_n)$ be one of its vertices chosen uniformly at random. Let $(\m_n,[\lab_n])\in\LMnsnun$ be the image of $(\q_n,v_n^\ooo)$ through the two-to-one mapping of Proposition~\ref{bijqvml} and let
$$\big(\s_n,(\f_n^e,\lab_n^e)_{e\in\vec F(\s_n)},(\bb_n^e)_{e\in\vec F(\s_n)}\big)$$
be the decomposition of $(\m_n,[\lab_n])$ appearing in Proposition~\ref{propdec}. We let $(\xi_n^e)_{e\in\vec E(\s_n)}$, $(m_n^e)_{e\in\vec F(\s_n)}$ and $(l_n^v)_{v\in V(\s_n)}$ be the three collections of integers from Proposition~\ref{propdec}. For all~$e$, we also denote by $(C_n^e,L_n^e)$ the contour pair of $(\f^e_n,\lab_n^e)$ as well as~$B_n^e$ the interpolation of~$\bb_n^e$. Finally, we define the rescaled versions of these objects
$$m_{(n)}^e \de \frac {2m_n^e + \xi_n^e}{2n},\qquad
\xi_{(n)}^e \de \frac{\xi_n^e}{\sqrt{2n}},\qquad
l_{(n)}^v \de \frac{l_n^v}{\g},$$
$$C_{(n)}^e \de \lp \frac {C_n^e(2ns)} {\sqrt{2n}} \rp_{0\le s \le m_{(n)}^e}\!,\ \ 
L_{(n)}^e \de \lp \frac {L_n^e(2ns)} {\g}\rp_{0\le s \le m_{(n)}^e}\!,\ \
B_{(n)}^e \de \lp \frac {B_n^e(\sqrt {2n}\, s)}{\g} \rp_{0\le s \le \xi_{(n)}^e}\! .$$
\nomenclature[6a0]{$C_{(n)}^e$, $L_{(n)}^e$, $B_{(n)}^e$}{rescaled coding functions corresponding to the half-edge~$e$}%
The goal of this section is to give the limit of the joint distribution of these processes. We first need to introduce the limiting processes.

\subsection{Brownian bridges, first-passage Brownian bridges, and Brownian snake}\label{secbbfpbbs}

We will work on the space $\K \de \bigcup_{x \in \R_+} \C([0, x],\R)$, endowed with the metric
$$d_\K (f,g) \de |\zeta(f) - \zeta(g) | + \sup_{y \ge 0} \lt f\big(y\wedge \zeta(f)\big)-g\big(y\wedge\zeta(g)\big)\rt,$$
where~$\zeta(f)$ denotes the only~$x$ such that $f \in \C([0, x],\R)$.

We call \emph{Brownian bridge} of length~$\xi$ from~$a$ to~$b$ a standard Brownian motion on $[0,\xi]$ started at~$a$, conditioned on being at~$b$ at time~$\xi$ (see for example \cite{billingsley68cpm, revuz99cma, bertoin03ptf, bettinelli10slr}). We also call \emph{first-passage Brownian bridge} of length~$m$ from~$a$ to~$b<a$ a standard Brownian motion on $[0,m]$ started at~$a$, and conditioned on hitting~$b$ for the first time at time~$m$. We refer the reader to~\cite{bettinelli10slr} for a proper definition of this conditioning, as well as for some convergence results of the discrete analogs.

The so-called \emph{Brownian snake's head} driven by a process $X\in \C([0, x],\R)$ may be defined as the process $(X(s),Z(s))_{0\le s \le x}$, where, conditionally given~$X$, the process~$Z$ is a centered Gaussian process with covariance function
\begin{equation*}
\cov\big(Z(s), Z(s')\big) = \inf_{s\wedge s'\le t\le s\vee s'} \big(X(t) - \underline X(t)\big).
\end{equation*}
We refer to \cite{legall99sbp, duquesne02rtl, bettinelli10slr} for more details about this process.

\subsection{Convergence of the encoding elements}

In the limit, we will see that only the schemes that maximize the cardinal of~$\vec E$ remain.
\begin{defi}\label{defscheme}
A scheme is \emph{dominant} if it has one vertex of degree exactly~$1$ and if all its other vertices are of degree exactly~$3$. Let~$\Sg_{p,1}^\star$ and~$\Sg_{p,2}^\star$ denote the sets of dominant schemes of type $(p,1)$ and $(p,2)$.
\nomenclature[5s31]{$\Sg_{p,q}^\star$}{set of dominant schemes of type $(p,q)$}%
\end{defi}

For instance, the two right-most schemes of Figure~\ref{sg11} are the two elements of $\Sg_{1,1}^\star$ in genus~$0$. Note that the degree $1$-vertex in the previous definition is necessarily an extremity of the root and that the dominant schemes are the ones whose number of edges is maximal (this is a simple consequence of Euler characteristic formula). Moreover, the condition on the degrees implies that every vertex of a dominant scheme is incident to at most one hole. In particular, the holes of a dominant scheme are well ``separated'' in the sense that their boundaries are pairwise disjoint simple loops, which are connected by some edges.

The compatibility condition on the previous collections of integers lead us to define, for a scheme~$\s$ with root~$e_*$, the set~$\mathcal{T}_\s$ of triples
$$\lp \lp m^e\rp_{e\in \vec F(\s)}, \lp\xi^e\rp_{e\in \vec E(\s)}, \lp l^v\rp_{v\in V(\s)} \rp \in \R_+^{\vec F(\s)} \times \R_+^{\vec E(\s)} \times \R^{V(\s)}$$
such that
\begin{itemize}
	\item $\sum_{e\in\vec F(\s)} m^e=1$,
	\item for all $e \in \vec E(\s)$, $\xi^{\bar e}=\xi^e$,
	\item for $0\le i \le p$, $\sum_{e\in\vec H_i(\s)} \xi^{e}=\si^i_\8$
	\item $l^{e_*^-}=0$.
\end{itemize}
We define the measure~$\mathcal{L}_\s$ on~$\mathcal{T}_\s$ as follows. For every $1\le i \le p$, we distinguish a half-edge $\vh_i\in\vec H_i(\s)$. We also consider an orientation~$\ori I(\s)$ of~$\vec I(\s)$, that is, a subset of~$\vec I(\s)$ containing exactly one half-edge among $\{e,\bar e\}$ for every $e\in\vec I(\s)$. There is then an obvious bijection between a part of
$$\R_+^{\vec F(\s)\bs\{e_*\}} \times \R_+^{\ori I(\s)\cup \bigcup_i \vec H_i(\s)\bs\{\vh_i\}} \times \R^{V(\s)\bs\{e_*^-\}}$$
and~$\mathcal{T}_\s$. The measure~$\mathcal{L}_\s$ is defined as the push-forward of the Lebesgue measure on this part.
%
%
We denote by~$p_a$ the density of a centered Gaussian variable with variance $a>0$, as well as~$-q_a$ its derivative:
$$p_a(x) \de \frac 1 {\sqrt{2\pi\, a}} \, \exp\lp - \frac {x^2} {2a} \rp \sand q_a(x) = \frac x a \, p_a \lp x \rp,\quad x\in\R.$$
We let~$\mu$ be the probability measure on $\bigcup_{\s\in\Sg_{p,1}^\star} \{\s\}\times \mathcal{T}_\s$ defined, for all measurable function~$\psi$, by
$$\mu(\psi) \de \frac1\Upsilon \sum_{\s\in\Sg_{p,1}^\star} \int_{\mathcal{T}_\s} d\mathcal{L}_\s \ 
		\psi\lp \s, \big( \lp m^e\rp, \lp\xi^e\rp, \lp l^v\rp \big) \rp
		\prod_{e\in{\vec F}(\s)} q_{m^e} \lp \xi^e \rp
		\!\!\!\!\prod_{e\in\ori{I}(\s)\cup\vec{B}(\s)}\!\!\!\! p_{(\kappa^e)^2\xi^e} \lp l^{e^+}\! - l^{e^-} \rp,$$
where
$$\kappa^e \de \begin{cases}
	1	&\text{ if } e \in \vec I(\s_\8)\\
	\sqrt 3	&\text{ if } e \in \vec B(\s_\8)
\end{cases}$$
and
\begin{equation*}\label{upsi}
\Upsilon \de \sum_{\s\in\Sg_{p,1}^\star} \int_{\mathcal{T}_\s} d\mathcal{L}_\s
		\prod_{e\in{\vec F}(\s)} q_{m^e} \lp \xi^e \rp
		\!\!\!\!\prod_{e\in\ori{I}(\s)\cup\vec{B}(\s)}\!\!\!\! p_{(\kappa^e)^2\xi^e} \lp l^{e^+}\! - l^{e^-} \rp
\end{equation*}
is a normalization constant. The measure~$\mu$ may also be described in terms of Gaussian Free Fields on the metric graphs corresponding to the schemes. We refer to the remark before Section~4.2 in~\cite{miermont09trm} for further details about this fact as it will not be used in this work. In fact, the precise expression of~$\mu$ will not be needed in what follows; we gave it for self-containment reasons and as it allows to compute the precise distribution of the random variable~$H$ appearing in Theorem~\ref{thmgeod}.

\begin{prop}\label{cvint}
The random vector
$$\lp \s_n, \big( m_{(n)}^e\big)_{e\in \vec F(\s_n)}, \big(\xi_{(n)}^e\big)_{e\in \vec E(\s_n)}, \big( l_{(n)}^v\big)_{v\in V(\s_n)}, \big( C_{(n)}^e, L_{(n)}^e \big)_{e\in \vec F(\s_n)}, \big( B_{(n)}^e \big)_{e\in \vec F(\s_n)} \rp$$
converges in distribution toward a random vector
$$\lp \s_\8, \lp m_\8^e\rp_{e\in \vec F(\s_\8)}, \lp\xi_\8^e\rp_{e\in \vec E(\s_\8)}, \lp l_\8^v\rp_{v\in V(\s_\8)}, \lp C_\8^e, L_\8^e \rp_{e\in \vec F(\s_\8)}, \lp B_\8^e \rp_{e\in \vec F(\s_\8)} \rp$$
whose law is defined as follows:
\begin{itemize}
	\item the vector $\lp \s_\8, \lp \lp m_\8^e\rp_{e\in \vec F(\s_\8)}, \lp\xi_\8^e\rp_{e\in \vec E(\s_\8)}, \lp l_\8^v\rp_{v\in V(\s_\8)}\rp \rp$ is distributed according to the probability measure~$\mu$;
	\item conditionally given this vector, 
	\begin{itemize}
		\item the processes $\lp C_\8^e, L_\8^e \rp$, ${e\in \vec F(\s_\8)}$ and $B_\8^e$, ${e\in \ori I(\s_\8)\cup \vec B(\s_\8)}$ are independent,
		\item the process $\lp C_\8^e, L_\8^e \rp$ has the law of a Brownian snake's head driven by a first-passage Brownian bridge of length~$m_\8^e$ from~$\xi_\8^e$ to~$0$,
		\item the process~$B_\8^e$ has the law of a Brownian bridge of length~$\xi_\8^e$ from~$l_\8^{e^-}$ to~$l_\8^{e^+}$, multiplied by the factor~$\kappa^e$,
		\item the bridges are linked through the relation $B_\8^{\bar e}(s) = B_\8^e(\xi_\8^e-s)$, $0\le s \le \xi_\8^e$, whenever $e \in \vec I(\s_\8)$.
	\end{itemize}
\end{itemize}
\end{prop}

The previous proposition is easily obtained by the method provided in~\cite{bettinelli10slr} (see in particular Proposition~7 and Section~5, as well as \cite[Proposition~7]{bettinelli11slr}); we leave the details to the reader.

The factor~$\kappa^e$ accounts for the fact that the steps of boundary bridges have a larger variance than the steps of interior bridges. This seemingly harmless factor causes some difficulties for the technical estimates of~\cite{bettinelli11slr}. Note also that this proposition is the reason why the factor $(8/9)^{1/4}$ appears. 

Let us emphasize at this stage that the limiting scheme~$\s_\8$ is a.s.\ dominant and, as such, possesses the properties observed after Definition~\ref{defscheme}.

We obtain a similar statement when performing the two-point mapping instead of the one-point mapping, up to a bias induced by the factor~$\lambda$ of Section~\ref{sectwopoint}. We will not need this statement in full details and will come back to it during Step~2.\ of the proof of Proposition~\ref{geodtyp}.

\begin{framed}
From now on, we apply Skorokhod's representation theorem and assume that the convergence of Proposition~\ref{cvint} holds almost surely. In particular, this entails that $\s_n=\s_\8$ for~$n$ large enough.
\end{framed}

\subsection{Proof of Theorem~\ref{cvthm}}

The whole proof of Theorem~\ref{cvthm} is obtained by adding the arguments of \cite{bettinelli10slr,bettinelli12tsl,bettinelli11slr}. The general strategy is borrowed from~\cite{legall07tss}, as well as from~\cite{miermont08sphericity} for the topology. We will not treat it here in full details and we refer the interested reader to these references for more information. We will only recall the main lines of reasoning and introduce the notions that will become useful in the next sections. We also exclude the case $(g,p)=(0,0)$, which is degenerate in some sense and requires special definitions. This is the original case treated in~\cite{legall07tss,miermont08sphericity}.

\subsubsection{Convergence of the metric space}\label{secconvms}
 
The first assertion is proved by the method of Le~Gall's pioneering paper~\cite{legall07tss}. Recall that~$\q_n$ is uniformly distributed over~$\Qnsn$, that~$v_n^\ooo$ is  uniformly distributed over~$V(\q_n)$ and that $(\m_n,[\lab_n])\in\LMnsnun$ is the labeled map corresponding to $(\q_n,v_n^\ooo)$. We arrange the corners of the internal face of~$\m_n$ according to the contour order, starting from the corner preceding the root: this gives a natural (non injective) ordering of the vertices of $\m_n$, which we write $\m_n(0)$, \ldots, $\m_n(2n+|\si_n|)$ with a slight abuse of notation. As the vertex set of~$\m_n$ corresponds to $V(\q_n)\backslash\{v_n^\ooo\}$, this also provides an ordering $\q_n(0)$, \ldots, $\q_n(2n+|\si_n|)$ of the vertices of~$\q_n$. (The fact that~$v_n^\ooo$ is left out will not be important when we take the scaling limit.) 

We may equip the set $[0,1]$ with the pseudo-metric $d_{(n)}$ defined by 
$$d_{(n)}(s,t)\de \lp\frac 9 {8n}\rp^{1/4} d_{\q_n}\Big(\q_n\big((2n+|\si_n|)\,s\big),\q_n\big((2n+|\si_n|)\,t\big)\Big),\qquad s,t \in \left\{0,\frac1{2n+|\si_n|},\ldots,1\right\},$$
and by linear interpolation between these values. It is easy to see that the metric space $[0,1]_{/\{d_{(n)}=0\}}$ equipped with the quotient metric is close in the Gromov--Hausdorff sense to the rescaled metric space corresponding to~$\q_n$. As the labels of~$\m_n$ represent the distances in~$\q_n$ to the point~$v_n^\ooo$, it turns out that it is possible to bound this pseudo-metric by an explicit function of $\lab_n$, which converges toward an explicit process, thanks to Proposition~\ref{cvint}. This entails that, up to extraction, $d_{(n)}$ converges toward a pseudo-metric $\disig$, the convergence holding jointly with the one of Proposition~\ref{cvint}. We define the equivalence relation associated with it by saying that $s \sim t$ if $\disig(s,t)=0$, and we set $\qis \de [0,1]_{/\sim}$. We also denote by $\qis(s)$ the equivalence class of~$s$ in~$\qis$. We can then show that this yields the first assertion of Theorem~\ref{cvthm}.
\nomenclature[6a4]{$\sim$, $\simeq$}{equivalence relation defining~$\qis$, equivalence relation defining~$\Mi$}%

\begin{framed}
From now on, we fix a subsequence $(n_k)_{k\ge 0}$ along which the previous convergence holds and we let $(\qis,\disig)$ be the corresponding limiting space. 
\nomenclature[6a1]{$(\qis,\disig)$}{Brownian surface}%
Recall also that, by Skorokhod's representation theorem, we assumed that the convergence of Proposition~\ref{cvint} held almost surely.
\end{framed}

\subsubsection{Seeing the limit as a quotient of a simpler ``map''}\label{secquot}

This is the longest and most difficult step. In some sense, we can view it as the continuous version of the mapping from Section~\ref{sec1pt}. The idea is to construct a continuous counterpart to~$\m_n$ and then to perform a continuous version of the mapping: because of the scaling, the counterpart of the arcs added in the discrete mapping is, in the continuous case, identifications of points.
 The main ingredients come from~\cite{legall07tss} but we need some extra arguments because of the more intricate combinatorial structure of the quadrangulations we consider. There are mainly two difficulties of different nature. The first difficulty arises from the more complex structure of the map~$\m_n$; we overcome this difficulty with the arguments of~\cite{bettinelli10slr,bettinelli12tsl}. The second technicality is caused by the boundary and is harder to grasp at first. The fact that the Brownian motions have different diffusion factors on the boundary edges and on the internal edges of the scheme (the factors~$\kappa^e$ from Proposition~\ref{cvint}) induces a technical complication; we use the arguments of~\cite{bettinelli11slr} to treat this problem.

We now properly define the continuous analog~$\Mi$ of~$\m_n$, which was briefly introduced during Section~\ref{secgeodbs}. Roughly speaking, we glue together Brownian forests coded by the $C_\8^e$'s according to the scheme structure. Precisely, let $e^{[1]}$, \ldots, $e^{[|\vec F(\s_\8)|]}$ be the half-edges of~$\vec F(\s_\8)$ arranged according to the contour order around~$f^\ooo$, starting from the root of~$\s_\8$. For every $s\in [0,1)$, there exists a unique $1 \le k \le {|\vec F(\s_\8)|}$ such that
$$\sum_{i=1}^{k-1} m^{e^{[i]}}_\8 \le s < \sum_{i=1}^{k} m^{e^{[i]}}_\8.$$
We let $e_s \de e^{[k]}$ and $\langle s \rangle \de s- \sum_{i=1}^{k-1} m_\8^{e^{[i]}} \in [0,m_\8^{e_s})$. By convention, we set $e_1=e^{[1]}$ and $\langle 1 \rangle =0$. We define the relation~$\simeq$ on $[0,1]$ as the coarsest equivalence relation for which $s \simeq t$ if one of the following occurs:
\begin{subequations}
\begin{flalign}
\hspace{\leftmargin}
&e_s=e_t\quad \text{ and }\quad  C_\8^{e_s}\langle s \rangle= C_\8^{e_s}\langle t \rangle=\inf_{[\langle s \rangle\wedge \langle t \rangle, \langle s \rangle\vee \langle t \rangle]}  C_\8^{e_s}\ ;&\label{idsamefor}\\
&e_s=\overline {e_t},\ \  C_\8^{e_s}\langle s \rangle = \underline { C_\8^{e_s}}\langle s \rangle,\ \  C_\8^{e_t}\langle t \rangle = \underline{C_\8^{e_t}}\langle t \rangle\ \text{ and }\ C_\8^{e_s}\langle s \rangle = \xi_\8^{e_t} - C_\8^{e_t}\langle t\rangle\ ;&\label{idfacefor}\\
&\langle s\rangle = \langle t\rangle = 0\quad\text{ and }\quad e_s^- = e_t^-&\label{idnode}
\end{flalign}
\end{subequations}
where we wrote $C_\8^{e_s}\langle s \rangle$ instead of $C_\8^{e_s}(\langle s \rangle)$ for short. Equation~\eqref{idsamefor} identifies numbers coding the same point in one of the forests, equation~\eqref{idfacefor} identifies the floors of forests ``facing each other'': the numbers~$s$ and~$t$ should code floor points (second and third equalities) of forests facing each other (first equality) and correspond to the same point (fourth equality). Finally, equation~\eqref{idnode} identifies the nodes. We set $\Mi \de [0,1]_{/\simeq}$ and, for $s$, $t \in [0,1]$, we let $\Mi(s)$ be the equivalence class of~$s$ in the quotient and $\Mi([s,t])\de\{\Mi(r)\, :\, r\in[s,t]\}$.
\nomenclature[6a3]{$\Mi$,$\fl$}{underlying continuous map, floor of~$\Mi$}%

\begin{figure}[ht!]
		\psfrag{a}[][][.8]{$a$}
		\psfrag{b}[][][.8]{$b$}
		\psfrag{B}[][][.8]{\textcolor{red}{$\lhb a, b \rhb$}}
		\psfrag{t}[][][.8]{$\tau_a$}
		\psfrag{i}[][][.8]{$h_1$}
		\psfrag{j}[][][.8]{$h_2$}
		\psfrag{f}[][]{$\fl$}
	\centering\includegraphics{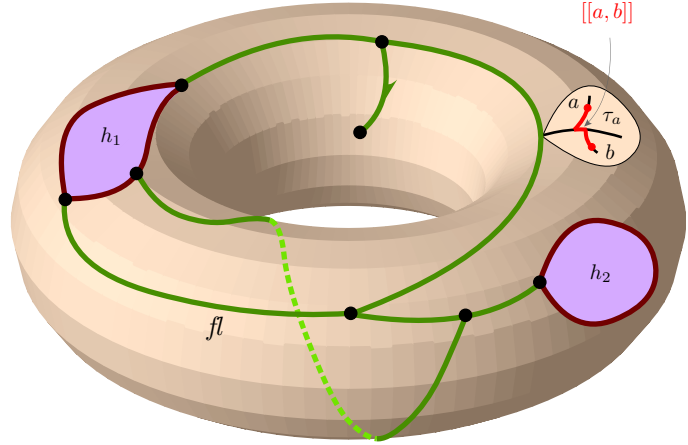}
	\caption{On this picture, we can see the floor $\fl$, an example of tree~$\tau_a$, and the set $\lhb a,b \rhb$.}
	\label{notm}
\end{figure}

The following notions will be used in what follows. They are the continuous counterpart of previous definitions. We call \emph{floor} of~$\Mi$ the set
$$\fl \de \Mi\lp \lb s \, :\, C_\8^{e_s}\langle s \rangle = \underline { C_\8^{e_s}}\langle s \rangle \rb \rp.$$
For $a=\Mi(s)\in\Mi\bs\fl$, let $l \de \inf\{t \le s \, :\, e_t=e_s,\ \underline { C_\8^{e_t}}\langle t \rangle=\underline { C_\8^{e_s}}\langle s \rangle \}$ and $r \de \sup\{t \ge s \, :\, e_t=e_s,\ \underline { C_\8^{e_t}}\langle t \rangle=\underline { C_\8^{e_s}}\langle s \rangle \}$. We set $\tau_a \de \Mi([l,r])$ and we call \emph{tree} of $\Mi$ a set of the form~$\tau_a$ for any $a\in\Mi\bs\fl$. We say that $[l,r]$ is the \emph{interval coding} the tree~$\tau_a$.  For $a$, $b$ in a tree of~$\Mi$, we denote by $\lhb a,b \rhb$ the range of the unique injective path linking~$a$ to~$b$.
\nomenclature[6a3b]{$\lhb a,b \rhb$}{range of the unique injective path linking~$a$ to~$b$}%
See Figure~\ref{notm}.

\begin{defi}\label{deforder}
The \emph{order} of a point $a\in\Mi$ is the integer
$$\big|\{s\in [0,1) \, : \, \Mi(s)=a\}\big|+ \un{\{a\in\fl, \exists s \, :\, a=\Mi(s),\, e_s\in\vec B(\s_\8)\}}.$$
A point of order~$1$ is called a \emph{leaf}.
\end{defi}

\begin{rem}
The indicator function in the previous definition accounts for the fact that the boundaries of the holes are not visited ``from the inside'' of the holes. This definition differs from the one we gave in Section~\ref{secgeodbs}; it is a simple exercise to verify that they are equivalent.
\end{rem}

Because the functions $C_\8^e$ are essentially Brownian motions, it is not hard to see that, a.s., all the points in~$\Mi$ are of order less than~$3$ and that, if~$U$ is uniformly distributed over~$[0,1]$, then $\Mi(U)$ is a.s.\ a leaf. Moreover, the floor contains no leaves.

The trivial fact that $\m_n((2n+|\si_n|)\,s)=\m_n((2n+|\si_n|)\,t) \Rightarrow d_{(n)}(s,t)=0$ easily passes to the limit in the sense that $s\simeq t \Rightarrow s\sim t$ (see \cite[Lemma~6]{bettinelli12tsl}). As a result, we may define a pseudo-metric and an equivalence relation on $\Mi$, which we still denote by~$\disig$ and~$\sim$, by setting $\disig\big(\Mi(s),\Mi(t)\big) \de \disig(s,t)$ and declaring $\Mi(s) \sim \Mi(t)$ if $s \sim t$. It allows us to see $\qis$ as a quotient of~$\Mi$: more precisely, the metric space $(\qis,\disig)$ is isometric to $\big({\Mi}_{/\sim},\disig\big)$.

\bigskip

We now describe how the points of~$\Mi$ are identified in the quotient $\Mi_{/\sim}$. We aim to assign labels to the points of~$\Mi$. We first define a labeling function on $[0,1]$ by setting
$$ \Lab(s)\de l_\8^{e_s^-}+ L_\8^{e_s}\langle s\rangle + B_\8^{e_s}\big( \xi_\8^{e_s} - \underline C_\8^{e_s}\langle s\rangle \big),\qquad 0\le s\le 1.$$
\nomenclature[6a5]{$\Lab$}{continuous label function, defined on $[0,1]$ or~$\Mi$}%
The reason for this definition is the observation we made at the end of Section~\ref{secrealvf}. Note that this function is continuous.

\begin{lem}\label{lemmin}
The set of points where $\Lab$ reaches its minimum is a.s.\ a singleton.
\end{lem}

We denote by~$s^\ooo$ the element of this singleton and we call \emph{base point} the corresponding point $\bp \de \qis(s^\ooo)$ in~$\qis$.
\nomenclature[6a6]{$s^\ooo$, $\rho^\ooo$}{$\operatorname{argmin} \Lab$, base point $\qis(s^\ooo)$ of~$\qis$}%
It is the continuous analog to the point~$v_n^\ooo$. (Recall that this point is linked to the vertices with minimum label in~$\m_n$.) As the labels in~$\m_n$ represent the distances to~$v_n^\ooo$ (up to an additive constant), a simple argument consisting in constructing a discrete approximation in~$\q_n$ of~$\qis(s)$ shows that 
\begin{equation}\label{disssooo}
\disig(s,s^\ooo) = \Lab(s) - \Lab(s^\ooo).
\end{equation}
The triangle inequality then entails that $s\sim t$ implies $\Lab(s)=\Lab(t)$. A fortiori, $s\simeq t$ implies $\Lab(s)=\Lab(t)$ so that we can view~$\Lab$ as a function on~$\Mi$ by setting $\Lab(\Mi(s))\de\Lab(s)$. The following proposition tells which points of~$\Mi$ are identified in~$\qis$. It is the continuous counterpart of the arc construction of Section~\ref{secrevmap}, roughly stating that points are identified if and only if they have the same label and all the points located inbetween when following the contour have greater labels.

\begin{prop}\label{propid}
Let $a=\Mi(s)$ and $b=\Mi(t)$ be two distinct points in~$\Mi$. Then $a\sim b$ if and only if~$a$ and~$b$ are leaves and
$$\Lab(a)=\Lab(b)\ge \Lab(c)$$
either for all $c\in\Mi([s\vee t,s\wedge t])$ or for all $c\in\Mi([s\wedge t,1]\cup[0,s\vee t])$.
\end{prop}

The two different sets appearing at the end come from the fact that there are two possible ways to follow the contour between~$a$ and~$b$: in some sense, one corresponds to linking~$a$ to~$b$ and the other one to linking~$b$ to~$a$. 

The proofs of these lemma and proposition are rather technical and use fine estimates on the Brownian snake. Happily, as the building blocks of our space are Brownian forests (the pieces coded by the sets $\overline{\{s<1\, :\, e_s=e\}}$), we may use the results of~\cite{bettinelli11slr} to quickly conclude.

\begin{pre}[Proof of Lemma~\ref{lemmin}]
Every half-edge $e\in\vec F(\s_\8)$ corresponds to a Brownian forest. They are not independent as the forests corresponding to~$e$ and~$\bar e$ when $e\in \vec I(\s_\8)$ share the same labels on their floors. However, me may consider these two forests as a single one having a density of trees multiplied by two. This merely changes the intensity of the Poisson process of the trees. With this convention, the edges of~$\s_\8$ correspond to independent forests so that they admit distinct minimums. The proof of \cite[Lemma~11]{bettinelli11slr} then allows to conclude that the minimum is reached only once a.s.
\end{pre}

\begin{pre}[Proof of Proposition~\ref{propid}]
The proof of this proposition is very similar to the proofs in the previously cited references. We will only give the main idea and refer the reader to these references for more details. The starting point is the following lemma. An \emph{increase point} for a function is a point~$t$ such that the function is greater than its value at~$t$ on a small interval of the form $[t-\eps,t]$ or $[t, t+\eps]$ for some $\eps>0$.
\begin{lem}\label{lemip}
The coding functions $s\in[0,1]\mapsto C_\8^{e_s}\langle s\rangle$ and~$\Lab$ do not share any increase points.
\end{lem}
This lemma is shown thanks to \cite[Lemma~12]{bettinelli11slr}, which states the result for a forest. Beware that, as explained in the latter reference, the result only holds for values of $\kappa^e\le\sqrt 3$, which is the case here. As the pair of functions considered here is a concatenation of pairs of functions for which the results holds, the result is straightforward. Some care is needed at the extremities of the intervals; we can show that such points are never increase points for the label process.

The second property we need is a property roughly stating the following. If~$a$ and~$b$ have the same label, if there exists a subtree~$\tau$ rooted at~$\rho$ such that $\inf_{\tau} \Lab < \Lab(a) < \Lab(\rho)$ and if, for infinitely many~$n$'s, there exists a geodesic from~$a_n$ to~$b_n$ that completely passes through~$\tau_n$, then $a\not\sim b$. Here, $a_n$, $b_n$ and~$\tau_n$ designate discrete approximations in~$\m_n$ of~$a$, $b$ and~$\tau$. See \cite[Lemma~15]{bettinelli11slr} for a rigorous statement. As explained in~\cite{bettinelli11slr}, in order to prove this property in our case, we only need an analog of \cite[Lemma~31]{bettinelli11slr}, which is easily obtained by decomposing~$\Lab$ according to the half-edges of~$\s_\8$.

These two properties together allow us to conclude as in \cite[Section~4.4]{bettinelli11slr}. The idea, coming from~\cite{legall07tss}, can be sketched as follows. The absence of common increase points entails with a little work that a point~$a$ of order at least~$2$ will be ``surrounded'' by subtrees arbitrarily close with a root having a label strictly greater than~$\Lab(a)$ and a minimum label strictly smaller than~$\Lab(a)$. As a result, it cannot be identified with any other point by the previous property. We then argue by contradiction and suppose that $a=\Mi(s)$ and $b=\Mi(t)$ are identified leaves but such that the property of the statement does not hold. After some work, it follows that on both $\Mi([s\vee t,s\wedge t])$ and $\Mi([s\wedge t,1]\cup[0,s\vee t])$ there exists a subtree with properties similar as above and we conclude that $a\not\sim b$, which is a contradiction.
\end{pre}

\subsubsection{Topology of the limit}\label{sectopo}

We follow the approach of~\cite{miermont08sphericity} to identify the topology of the limit. The first observation is that the metric space associated with~$\q_n$ is not far from being homeomorphic to~$\Sgp$. Informally, we fill the internal faces with small cells and the holes with annuli. The rigorous way to proceed is to consider, for every internal face a copy of the space $[0,1]^3\bs ((0,1)^2 \times [0,1))$ and for every hole~$h_i$ a space $(\mathcal G_i \times [0,1]) \bs (\mathring{\mathcal G_i} \times [0,1])$ where $\calG_i$ is a regular $2\sigma_n^i$-sided polygon with unit length edges embedded in~$\R^2$. We endow these spaces with the intrinsic metric inherited from the Euclidean metric, rescaled by $(9/8n)^{1/4}$. We can then define a surface (with a boundary) by identifying, according to the map structure, the boundaries of these spaces and endow it with the quotient metric. The metric space we obtain is homeomorphic to~$\Sgp$ and the choice of cells and annuli we made ensures that, on the set of points corresponding to the vertices of the map, the metric coincides with the graph metric of the map (in other words, we did not create any shortcuts as it would have been the case if we simply had glued Euclidean polygons, for instance). As a result, this surface is at Gromov--Hausdorff distance $o(1)$ from $( V(\q_{n}),(9/8n)^{1/4} d_{\q_{n}})$. See \cite[Section~5.2]{bettinelli11slr} for more details.

The second step is a criterion ensuring that the Gromov--Hausdorff limit of a sequence of metric spaces all homeomorphic to some surface is also homeomorphic to this surface. This criterion is called regularity and was introduced by Whyburn~\cite{whyburn35sls} and later studied by Begle~\cite{begle44rc}. We will not go into much detail as the arguments are exactly the same as in~\cite{bettinelli12tsl,bettinelli11slr}. The only difference is that, a priori, the regularity criterion only works for~$\Sgpof g 0$ (\cite[Proposition~19]{bettinelli12tsl}) and for~$\Sgpof{0}{1}$ (\cite[Proposition~16]{bettinelli11slr}). We may adapt it to work for~$\Sgp$ by noticing that~$\Sgp$ can be constructed by removing from~$\Sgpof{g}{0}$ the interior of~$p$ disjoint subspaces homeomorphic to~$\Sgpof{0}{1}$. We first ``fill'' the holes (formally by changing $(\mathcal G_i \times [0,1]) \bs (\mathring{\mathcal G_i} \times [0,1])$ into $(\mathcal G_i \times [0,1]) \bs (\mathring{\mathcal G_i} \times [0,1))$) and show a convergence toward a surface homeomorphic to~$\Sgpof{g}{0}$ thanks to the first theorem. We then show that the~$p$ closed holes (corresponding in the filled surface to $\mathcal G_i \times \{1\}$) converge to~$p$ disjoint disks thanks to the second theorem. In a few words, this comes from the fact that the holes of~$\q_n$ correspond to the holes of~$\m_n$ and, if they were not converging to disjoint disks, in the limit, we would either have holes ``sharing'' a vertex in the scheme or we would have points of the floor of~$\Mi$ identified together. Neither of these situations can happen: as observed after Definition~\ref{defscheme}, there are no vertices incident to two different holes in a dominant scheme, and Proposition~\ref{propid} prohibits the latter situation.

\subsubsection{Hausdorff dimension of the limit}

To compute the Hausdorff dimension of the limit, we use the same method as in~\cite{legall09scaling}. The idea is roughly the following. To prove that the Hausdorff dimension is less than~$4$, we use the fact that the processes $L_\8^e$ of Proposition~\ref{cvint} are almost surely $\alpha$-H\"older for all $\alpha \in (0 , 1/4)$, yielding that the canonical projection from $([0,1],|\cdot|)$ to $(\qis,\disig)$ is also $\alpha$-H\"older for the same values of~$\alpha$. To prove that it is greater than~$4$, we show that the size of the balls of diameter~$\delta$ is of order~$\delta^4$. To see this, we first bound from below the distances in terms of label variation along the branches of~$\Mi$, and then use twice the law of the iterated logarithm: this tells us that, for a fixed $s\in [0,1]$, the points outside of the set $[s-\delta^4,s+\delta^4]$ code points that are at distance at least~$\delta^2$ from~$\qis(s)$ in~$\Mi$, so that their distance from~$\qis(s)$ is at least~$\delta$ in the map. See \cite[Section~6.4]{bettinelli10slr} for a complete proof.

We now turn to the boundary. Let $\pii:\Mi\to\qis$ be the canonical projection.
\nomenclature[6a4b]{$\pii$}{canonical projection $\Mi\to\qis$}%
A straightforward adaptation of \cite[Proposition~21]{bettinelli11slr} shows that the boundary of~$\qis$ is
\begin{equation}\label{boundary}
\partial\qis = \pii\lp\Mi\lp{\{s\, :\, e_s\in\vec B(\s_\8)\}}\rp \cap \fl\rp.
\end{equation}
\nomenclature[6a2]{$\partial\qis$}{boundary of $(\qis,\disig)$}%
For every connected component of the boundary, we may use the same method as in~\cite{bettinelli11slr} to conclude that its dimension is a.s.\ equal to~$2$.

\section{Typical geodesics from the base point}\label{sectypgeod}

From now on, we work on the set of full probability where~$\Lab$ reaches its minimum only once. Recall that~$\pii$ denotes the canonical projection from~$\Mi$ to~$\qis$, that~$s^\ooo$ is the unique point where~$\Lab$ reaches its minimum and that we set $\bp \de \qis(s^\ooo)$. We also include back the case $(g,p)=(0,0)$ from now on. In this case, $(\m_n,[\lab_n])\in\mM_{n,\emptyset}^1$ is a uniform labeled tree and its rescaled contour pair converges toward the Brownian snake $(\ee,Z)$ driven by a normalized Brownian excursion. Thus~$\Mi$ is Aldous's CRT coded by the excursion~$\ee$ and $\Lab\de Z$. Lemma~\ref{lemmin} and Proposition~\ref{propid} also hold (see~\cite{legall07tss}).

\subsection{Uniqueness}\label{secuniq}

\begin{prop}\label{geodtyp}
Let $S$ be uniformly distributed over $[0,1]$ and independent of $(\qis,\disig)$. Then, almost surely, there is only one geodesic from the base point~$\bp$ to~$X\de\qis(S)$ in~$\qis$. 
\end{prop}

\begin{pre}
This was shown in~\cite{miermont09trm} in the slightly different context of Boltzmann maps. We reformulate the ideas to our setting.

It will be sufficient to show that, for any rational numbers $0<\alpha<\beta$, the event $\EE_{\alpha,\beta}$ that there exist, for every $u\in[\alpha,\beta]$, at least two different points~$y$ and $y'\in \qis$ such that $\disig(\bp,y)=\disig(\bp,y')=u$ and $\disig(\bp,y)+ \disig(y,X)=\disig(\bp,y')+ \disig(y',X)=\disig(\bp,X)$ has probability~$0$. For the remaining of the proof, we fix rational numbers $0<\alpha<\beta$ and we consider a random variable~$U$ uniformly distributed over $[\alpha,\beta]$ and independent of all other variables.

The idea is that the previous points~$y$ and~$y'$ correspond to global minimums of the labels on the interface between the two internal faces of the map obtained by a continuous analog of the two-point mapping of Section~\ref{sectwopoint}. As the labels are essentially Brownian, this almost surely never happens. Let us proceed rigorously.

\paragraph{Step 1.}
We go back to the discrete picture and start with the construction of a suitable approximation of~$X$ in~$\q_n$. We want to find a sequence $(S_n)$ of integers such that $\q_n(S_n)$ is roughly uniformly distributed over~$V(\q_n)$ and $S_n/2n\to S$, at least along a subsequence. When~$\m_n$ is a tree, this is quite easy. In the general case, it is a bit trickier; we proceed as follows. Let~$w$ be an arbitrarily chosen point at $d_{\q_n}$-distance~$1$ from~$v_n^\ooo$. We claim that we can find a sequence~$(S_n)$ of integers such that $X_n\de\q_n(S_n)$ is uniformly distributed over~$V(\q_n)\bs\{v_n^\ooo,w\}$ and $S_n/2n\to S$ along some subsequence. To do so, we will associate with every of these $n+|\sigma_n|-p-2g$ vertices exactly two unit length left-open subintervals of $[0,2(n+|\sigma_n|-p-2g)]$ as follows. It is possible to remove $2g+p$ edges from~$\m_n$ without disconnecting the map; we obtain a plane tree. In the resulting tree, we bijectively associate with every vertex different from~$w$ an edge in such a way that the edge is incident to the corresponding vertex. Then, in the contour of~$f^\ooo$, an edge of~$\m_n$ corresponds to one or two unit length subintervals of $[0,2n+|\sigma_n|]$ (one if the edge is incident to a hole, two otherwise). We suppose that $|\sigma_n|\ge p+2g$, which happens for large~$n$, and we associate with every vertex of~$\m_n$ the corresponding subintervals. We associate with the vertices corresponding to only one subinterval a second subinterval arbitrarily chosen among the remaining ones. Let~$k$ be such that $2(n+|\sigma_n|-p-2g)S\in (k,k+1]$. If the vertex corresponding to $(k,k+1]$ is either~$\m_n(k)$ or~$\m_n(k+1)$, we set~$S_n$ to~$k$ or~$k+1$ accordingly. Otherwise, we arbitrarily set~$S_n$ such that~$\m_n(S_n)$ corresponds to $(k,k+1]$. With this construction, $X_n$ is clearly uniformly distributed over~$V(\q_n)\bs\{v_n^\ooo,w\}$ as claimed. Moreover, the probability that we cannot extract a subsequence along which $S_n/2n\to S$ is bounded by
$$\Pb\lp\liminf_{n\to\8} \big\{|S_n-2(n+|\sigma_n|-p-2g)S|>1\big\}\rp\le \liminf_{n\to\8} \frac{|\sigma_n|}{2(n+|\sigma_n|-p-2g)}=0.$$
So, up to discarding a zero-probability event, we may suppose that $S_n/2n\to S$ along a subsequence of $(n_k)_{k\ge 0}$.

\paragraph{Step 2.}
We look at the scaling limit of the map obtained by the two-point mapping. We set $\lambda_n\de \lfloor (8n/9)^{1/4}\, U \rfloor$. If $1\le\lambda_n\le d_{\q_n}(v_n^\ooo,X_n)-1$, it is possible to define $(\m^{\ooo\ooo}_n,[\lab^{\ooo\ooo}_n])\in\LMnsnde$ as the labeled map corresponding to $(\q_n,(v_n^\ooo,X_n), \lambda_n )$ via the two-point mapping of Section~\ref{sectwopoint}.

We decompose a labeled map $(\m,[\lab])\in \LMnsnde$ into a scheme and a collection of integers and continuous functions as in Section~\ref{secsl}. We denote by $\mathcal V_{(n)}(\m,[\lab])$ the corresponding rescaled vector (the analog of the one from Proposition~\ref{cvint}). By convention, we set $\mathcal V_{(n)}(\m^{\ooo\ooo}_n,[\lab^{\ooo\ooo}_n])\de 0$ whenever $(\m^{\ooo\ooo}_n,[\lab^{\ooo\ooo}_n])$ is not defined (that is, whenever $\lambda_n \notin\{1,\ldots,d_{\q_n}(v_n^\ooo,X_n)-1\}$). We compare the distribution of $\mathcal V_{(n)}(\m^{\ooo\ooo}_n,[\lab^{\ooo\ooo}_n])$ with the distribution of $\mathcal V_{(n)}(M_n,[L_n])$, where $(M_n,[L_n])$ denotes a random variable uniformly distributed over $\LMnsnde$: For any bounded measurable function~$F$ such that $F(0)=0$, a straightforward computation yields that
$$\frac{\mathbb E\big[F\big(\mathcal V_{(n)}(\m^{\ooo\ooo}_n,[\lab^{\ooo\ooo}_n])\big) \big]}{\Pb\big(\lambda_n<d_{\q_n}(v_n^\ooo,X_n)\big)}=
\frac{\mathbb E\big[F\big(\mathcal V_{(n)}(M_n,[L_n])\big)\,\Delta_n\big(\mathcal V_{(n)}(M_n,[L_n])\big)\big]}{\E{\Delta_n\big(\mathcal V_{(n)}(M_n,[L_n])\big)}},$$
where $\Delta_n$ is defined as follows. For $(\m,[\lab])\in \LMnsnde$, the number $\Delta_n(\mathcal V_{(n)}(\m,[\lab]))$ is $(8n/9)^{1/4}$ times the probability that~$\lambda_n$ equals the~$\lambda$ corresponding to $(\m,[\lab])$ via the two-point mapping. It is a simple piecewise linear function of $(9/8n)^{1/4}\lambda$, equal to~$0$ on $(-\8,\alpha-(9/8n)^{1/4}]\cup[\beta,\8)$ and equal to $1/(\beta-\alpha)$ on $[\alpha,\beta-(9/8n)^{1/4}]$. Moreover, by Lemma~\ref{lemlambda}~($iii$), $(9/8n)^{1/4}(\lambda-1)$ may be written as a deterministic continuous function~$\Lambda$ of $\mathcal V_{(n)}(\m,[\lab])$. The explicit expression of~$\Lambda$ is quite intricate and will not be needed: loosely speaking, $\Lambda$ assigns to a vector the minimum of the labels on the interface between the two faces minus the minimum of the labels inside~$f^\ooo$.

Using the same machinery as in Section~\ref{secsl}, we obtain that the random vector $\mathcal V_{(n)}(M_n,[L_n])$ converges to a random vector~$\mathcal V_\8$ whose precise distribution will not be needed. We easily deduce that $\mathcal V_{(n)}(\m^{\ooo\ooo}_n,[\lab^{\ooo\ooo}_n])$ converges weakly toward a random vector~$\mathcal V_\8^{\ooo\ooo}$ with distribution
$$\Pb\big(U>\disig(\bp,X)\big)\, \delta_0 +\Pb\big(U<\disig(\bp,X)\big)\, \lL,$$
where~$\lL$ is the probability distribution given by
\begin{equation}\label{bias}
\lL(G) = \E{G(\mathcal V_\8)\,\frac{\un{\{\alpha<\Lambda(\mathcal V_\8) <\beta\}}}{\Pb(\alpha<\Lambda(\mathcal V_\8) <\beta)}}
\end{equation}
for any measurable function~$G$.
In particular, the distribution~$\lL$ is absolutely continuous with respect to the distribution of~$\mathcal V_\8$ so that every property holding a.s.\ for~$\mathcal V_\8$ also holds a.s.\ under~$\lL$.

\paragraph{Step 3.}
We now work on the event $\EE_{\alpha,\beta}$. We denote by~$y$ and~$y'$ two different points satisfying the relations appearing in the definition of $\EE_{\alpha,\beta}$ for $u=U$ and we let~$t$ and~$t'$ be such that $y=\qi(t)$ and $y'=\qi(t')$. We construct discrete approximations of~$y$ and~$y'$ as follows.

As $0<U<\disig(\bp,X)$, we see that $(\m^{\ooo\ooo}_n,[\lab^{\ooo\ooo}_n])$ is always defined for~$n$ large enough.
We denote by~$f_n^\ooo$ and~$f_n^{\ooo\ooo}$ its two internal faces corresponding respectively to~$v_n^\ooo$ and~$X_n$. We want to construct an approximation~$y_n$ on the interface $f_n^\ooo \cap f_n^{\ooo\ooo}$. Let $z_n\de \q_n(\lfloor 2nt \rfloor)$. If~$z_n\in f_n^\ooo \cap f_n^{\ooo\ooo}$, then we set~$y_n\de z_n$. Otherwise, let suppose for now that $z_n\in f_n^{\ooo\ooo}$. We consider a geodesic from~$z_n$ to~$v_n^\ooo$ and we let~$y_n$ be a point on this geodesic that belongs both to~$f_n^\ooo$ and~$f_n^{\ooo\ooo}$. By construction and by Lemma~\ref{lemlambda}~($ii$), we have $\lambda_n\le d_{\q_n}(v_n^\ooo,y_n) \le d_{\q_n}(v_n^\ooo,z_n)$. The triangle inequality yields that $d_{\q_n}(y_n,z_n) = o(n^{1/4})$ and another application of Lemma~\ref{lemlambda} gives
\begin{equation}\label{lnyn}
0\le \lab^{\ooo\ooo}_n(y_n) - \min_{f_n^\ooo\cap f_n^{\ooo\ooo}} \lab^{\ooo\ooo}_n = o(n^{1/4}).
\end{equation}
Finally, if~$z_n\in f_n^\ooo$, we construct~$y_n$ in a symmetrical way and obtain the same equations by similar manipulations.

We also construct the points~$z_n'$ and~$y_n'$ by applying the same method to~$t'$ instead of~$t$. As $(9/8n)^{1/4}\, d_{\q_n}(z_n,z_n') \to \disig(y,y')$, we also obtain that $(9/8n)^{1/4}\,d_{\q_n}(y_n,y_n')\to \disig(y,y')>0$. 

The interface $f_n^\ooo\cap f_n^{\ooo\ooo}$ is a ``cycle'' in the map~$\m_n^{\ooo\ooo}$. We parametrize this cycle by $[0,1]$, the point~$0$ being for example a node chosen in a deterministic way from the scheme. We denote by $r_{(n)}\in [0,1)$ and $r'_{(n)}\in [0,1)$ the real numbers corresponding to~$y_n$ and~$y'_n$. Up to further extraction, we may suppose that $r_{(n)}\to r$ and $r_{(n)}'\to r'$. As $(9/8n)^{1/4}\,d_{\q_n}(y_n,y_n')$ converges toward a positive number, we have $r\neq r'$: This is obtained by the argument we used in Section~\ref{secquot}. In short, we may construct an equivalence relation from the limiting vector~$\mathcal V_\8^{\ooo\ooo}$ by a construction similar to the one of Section~\ref{secquot} and show that two points identified by this relation are also identified in~$\qis$. In other words, because two corners identified in~$\m_n^{\ooo\ooo}$ correspond to the same point in~$\q_n$, the Brownian surface~$\qis$ may be seen as a quotient of a continuous analog of~$\m_n^{\ooo\ooo}$. This is very similar to what we did in Section~\ref{secquot} so that we leave the details to the reader.

Now, Step~$2$ implies that the labels on the interface $f_n^\ooo\cap f_n^{\ooo\ooo}$ converge in distribution, once properly rescaled, toward a concatenation of some functions of~$\mathcal V_\8^{\ooo\ooo}$ (up to an irrelevant time scaling because of the parametrization by $[0,1]$). Let us denote by $B_\8^{\ooo\ooo}:[0,1]\to\R$ this limiting function. Equation~\eqref{lnyn} yields that $B_\8^{\ooo\ooo}(r)=B_\8^{\ooo\ooo}(r')=\min B_\8^{\ooo\ooo}$. As, on the event $\{\mathcal V_\8^{\ooo\ooo} \neq 0\}$, the process~$B_\8^{\ooo\ooo}$ is, up to the bias~\eqref{bias}, a Brownian bridge, this means that we are on an event of probability zero, as wanted.
\end{pre}

\subsection{Simple geodesics}\label{secsg}

We use in this section the terminology and ideas of~\cite{legall08glp}, developed in the spherical case and readily adaptable to our more general setting. We define simple geodesics and give some elementary properties. Recall that $\disig(s^\ooo,s) = \Lab(s) - \Lab(s^\ooo)$. For $s$, $t\in[0,1]$, we set
$$\st s t \de \begin{cases} [s,t]	&\text{ if } s \le t, \\ [s,1] \cup [0,t]	& \text{ if } t < s. \end{cases}$$
\nomenclature[7a1]{$\st s t$}{$[s,t]$ if $s \le t$ or $[s,1] \cup [0,t]$ if $t < s$}%

\begin{defi}
The \emph{simple geodesic} of index~$s\in [0,1]$ is the path~$\Phi_s$ defined by
$$\Phi_s(w)\de \qis\lp \sup\lb r\, :\, \inf_{\st r s} \Lab = \Lab(s^\ooo)+ w\rb \rp,\qquad 0\le w\le \disig(s^\ooo,s).$$
\end{defi}
\nomenclature[7a2]{$\Phi_s$}{simple geodesic from~$\bp$ of index~$s$}%

It is not very hard to see that~$\Phi_s$ is a geodesic from~$\rho^\ooo$ to~$\qis(s)$ (see \cite[Section~4.1]{legall08glp} for more details). If we set $t\de \sup\{ r\, :\, \inf_{\st r s} \Lab = \Lab(s^\ooo)+ w\}$ then we must have $\Lab(t) = \inf_{\st {t} s} \Lab$ by continuity. If $w< \disig(s^\ooo,s)$, then $t\neq s$ so that~$t$ is an increase point of~$\Lab$ and, by Lemma~\ref{lemip}, $\Mi(t)$ is a leaf. This entails that~$\Phi_s$ only visit images of leaves, except possibly at its end~$\qis(s)$.
Moreover, if $t'\de \inf\{ r\, :\, \inf_{\st s r} \Lab = \Lab(s^\ooo)+ w\}$ then Proposition~\ref{propid} implies that $t\sim t'$. This provides a dual definition of~$\Phi_s$:
$$\Phi_s(w)=\qis\lp \inf\lb r\, :\, \inf_{\st s r} \Lab = \Lab(s^\ooo)+ w\rb \rp,\qquad 0\le w\le \disig(s^\ooo,s).$$

Observe that, if~$t''$ is such that $\Lab(t'') = \inf_{\st {t''} s} \Lab$, then $t''\sim \sup\{ r\, :\, \inf_{\st r s} \Lab = \Lab(s^\ooo)+ (\Lab(t'')- \Lab(s^\ooo))\}$. As a result, the range of~$\Phi_s$ is given by
$$\Phi_s\big([0,\disig(s^\ooo,s)]\big)=\qis\lp \lb t \, : \, \Lab(t) = \inf_{\st t s} \Lab \rb\rp=\qis\lp \lb t \, : \, \Lab(t) = \inf_{\st s t} \Lab \rb\rp.$$

Finally, let $s\sim s'$ be two distinct numbers. If $s\not\simeq s'$ then, by Proposition~\ref{propid}, we have either $\Lab(s)=\Lab(s')=\inf_{\st s {s'}}\Lab$ or $\Lab(s)=\Lab(s')=\inf_{\st {s'} s}\Lab$; thus $\Phi_s=\Phi_{s'}$. If $s\simeq s'$ then~$s$ and~$s'$ cannot be increase points of~$\Lab$, so that both $\inf_{\st s {s'}}\Lab<\Lab(s)$ and $\inf_{\st {s'} s}\Lab<\Lab(s)$. It is easy to see that $\Phi_s(w)=\Phi_{s'}(w)$ if and only if $w\in [0, \inf_{\st s {s'}}\Lab\vee \inf_{\st {s'} s}\Lab-\Lab(s^\ooo)]\cup\{\Lab(s)-\Lab(s^\ooo)\}$. In particular, $\Phi_s\neq\Phi_{s'}$.

\bigskip

Let $(s_i)_{i \ge 0}$ be a sequence of i.i.d.\ random variables uniformly distributed over~$[0,1]$. We denote by $a_i \de \Mi(s_i)$ and $x_i \de \qis(s_i)$ the corresponding points in~$\Mi$ and~$\qis$. By Proposition~\ref{geodtyp}, it is immediate that almost surely, for all~$i$, $\Phi_{s_i}$ is the only geodesic from~$\bp$ to~$x_i$. Note also that, almost surely, the set $\{s_i\, : \, i\ge 0\}$ is a dense subset of $[0,1]$. Up to discarding yet another event of zero probability, we will suppose that both previous properties hold. 

In what follows, we will show that it is actually sufficient to know that the geodesics to the~$x_i$'s are simple geodesics in order to conclude that all the geodesics from~$\bp$ are simple geodesics.

\section{General geodesics from the base point}\label{secgengeod}

\subsection{Geodesics from the base point to images of leaves}\label{secgeodleaf}

In this section, we treat the case of the image under~$\pii$ of the set of leaves in~$\Mi$. Let $s \in [0,1]$ be such that $a\de \Mi(s)$ is a leaf. We will show that~$\Phi_s$ is the only geodesic from~$\bp$ to $x\de \qis(s)$ in~$\qis$. There is nothing to prove in the case $s=s^\ooo$, so that we suppose $s \neq s^\ooo$. The first step is to find points~$a_i$ and~$a_j$ before and after~$a$ (in the contour of~$\Mi$), arbitrarily close to~$a$ and such that no geodesics from~$\bp$ to~$x$ intersect the set $\pii(\lhb a_i,a_j \rhb)$. Note that, as~$a$ is a leaf, the tree~$\tau_a$ is well defined.

\begin{lem}\label{passpas}
Let $\eps >0$ be such that $[s-\eps,s+\eps]$ is included in the interval coding~$\tau_a$. There exist~$i$ and~$j$ satisfying $s-\eps < s_i < s < s_j < s+\eps$ and, for all $b \in \lhb a_i,a_j \rhb$,
$$\disig(\bp,a) < \disig(\bp,b) + \disig(b,a).$$
\end{lem}

\begin{figure}[ht!]
		\psfrag{a}[][]{$a$}
		\psfrag{b}[][]{$b_\eta$}
		\psfrag{c}[][]{$c_\xi$}
		\psfrag{i}[][]{$a_i$}
		\psfrag{j}[][]{$a_j$}
		\psfrag{l}[][]{$l_\xi$}
		\psfrag{r}[][]{$r_\xi$}
		\psfrag{u}[][]{$s$}
		\psfrag{s}[][]{$s_i$}
		\psfrag{t}[][]{$s_j$}
		\psfrag{e}[][]{$t_\eta$}
		\psfrag{A}[][]{\textcolor{red}{$\lhb a_i,a_j \rhb$}}
		\psfrag{x}[][]{$\xi$}
	\centering\includegraphics{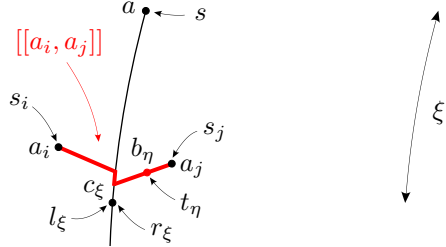}
	\caption{Notation used in the proof of Lemma~\ref{passpas}.}
	\label{aaiaj}
\end{figure}

\begin{pre}
We argue by contradiction and suppose that there exists an $\eps >0$ such that $[s-\eps,s+\eps]$ is included in the interval coding~$\tau_a$, and that for all~$i$ and~$j$ satisfying $s-\eps < s_i < s < s_j < s+\eps$, we can find $b\in \lhb a_i, a_j \rhb$ for which $\disig(\bp,a) = \disig(\bp,b) + \disig(b,a)$. For $0\le\xi\le C_\8^{e_s}\langle s \rangle - \underline { C_\8^{e_s}}\langle s \rangle$, we define
$$l_\xi \de \sup\{t \le s \, : \, C_\8^{e_t}\langle t \rangle = C_\8^{e_s}\langle s \rangle -\xi \} \sand r_\xi \de \inf\{t \ge s \, : \, C_\8^{e_t}\langle t \rangle = C_\8^{e_s}\langle s \rangle -\xi \}$$
and we set $c_\xi \de \Mi(l_{\xi})=\Mi(r_{\xi})$. See Figure~\ref{aaiaj}. As~$a$ is a leaf, there exists~$\xi_0>0$ such that $s-\eps < l_\xi < s < r_\xi < s+\eps$ as soon as $\xi\in (0,\xi_0]$. Until further notice, we fix such a~$\xi$.

Let $0 < \eta < (s-l_\xi) \wedge (r_\xi - s)$ and $i$, $j$ be such that $l_\xi \le s_i < l_\xi + \eta$ and $r_\xi -\eta < s_j \le r_\xi$. By hypothesis, we can find $b_\eta \in \lhb a_i, a_j \rhb$ such that $\disig(\bp,a) = \disig(\bp,b_\eta) + \disig(b_\eta,a)$. We can furthermore find $t_\eta \in [l_\xi, s_i] \cup [s_j, r_\xi]$ such that $b_\eta \de \Mi(t_\eta)$. As $\eta \to 0$, $t_\eta$ converges, up to extraction, either to~$l_\xi$ or to~$r_\xi$. This implies that $\disig(c_\xi,a) = \disig(\bp,a) - \disig(\bp,c_\xi)=\Lab(a)-\Lab(c_\xi)$, by~\eqref{disssooo}. Moreover, the classical ``cactus bound'' (see~\cite{curienlgm13cactusI} or \cite[Lemma~20]{bettinelli10slr} for a similar statement) states that $\disig(c_\xi,a) \ge \Lab(c_\xi) + \Lab(a) - 2\min_{\lhb c_\xi,a \rhb} \Lab$, so that $\Lab(c_\xi)= \min_{\lhb c_\xi,a \rhb} \Lab$. As a result, the process $\xi \in [0, \xi_0] \mapsto \Lab(c_\xi)$ is non-increasing. Basic properties of Brownian snakes show that this is a contradiction.
\end{pre}

We will need to consider sets of the form $\Mi([s_i,s_j])$. The following lemma identifies the topological boundary of the image in~$\qis$ of such a set.

\begin{lem}\label{boundt1t2}
Let~$\tau$ be a tree of~$\Mi$ coded by $[l,r]$ and $l \le t_1 \le t_2 \le r$. We set $w' \de \inf \{w\, : \, \Phi_{t_1}(w) \neq \Phi_{t_2}(w)\}$. Then the topological boundary in~$\qis$ of the subset $\qis([t_1,t_2])$ is
$$\partial\big( \qis([t_1,t_2])\big)=\pii\big(\lhb \Mi(t_1), \Mi(t_2) \rhb \big) \cup \Phi_{t_1}\big([w',\disig(s^\ooo,t_1)]\big) \cup \Phi_{t_2}\big([w',\disig(s^\ooo,t_2)]\big).$$
\end{lem}

\begin{figure}[ht!]
		\psfrag{a}[][]{$a$}
		\psfrag{b}[][]{$b_\eta$}
		\psfrag{c}[][]{$c_\xi$}
		\psfrag{i}[][]{$a_i$}
		\psfrag{j}[][]{$a_j$}
		\psfrag{l}[][]{$l_\xi$}
		\psfrag{r}[][]{$r_\xi$}
		\psfrag{u}[][]{$s$}
		\psfrag{s}[][]{$t_1$}
		\psfrag{t}[][]{$t_2$}
		\psfrag{e}[][]{$t_\eta$}
		\psfrag{A}[r][r]{\textcolor{red}{$\pi([[\Mi(t_1), \Mi(t_2)]])$}}
		\psfrag{1}[][]{\textcolor{blue}{$\Phi_{t_1}$}}
		\psfrag{2}[][]{\textcolor{blue}{$\Phi_{t_2}$}}
		\psfrag{x}[][]{$\xi$}
	\centering\includegraphics[width=55mm]{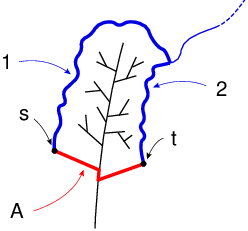}
	\caption{Topological boundary of $\qis([t_1,t_2])$.}
	\label{aaiajb5}
\end{figure}

In other words, the boundary of $\qis([t_1,t_2])$ is composed of three parts: the image of the set $\lhb \Mi(t_1), \Mi(t_2) \rhb$ and the ranges of~$\Phi_{t_1}$ and~$\Phi_{t_2}$ after they split. See Figure~\ref{aaiajb5}.

\begin{pre}
The properties stated during Section~\ref{secsg} show that
$$\Phi_{t_1}\big([w',\disig(s^\ooo,t_1)]\big)=\qis\lp \Big\{ t \in [t_1,t_2] \, : \, \Lab(t) = \inf_{[t_1,t]} \Lab \Big\} \rp.$$
Let~$t \in [t_1,t_2]$ be such that $\Lab(t) = \inf_{[t_1,t]} \Lab$, and let us define $t'\de \sup\{r\, : \, \inf_{\st r {t_1}}\Lab = \Lab(t)\}$. Then $t' \eqq t$, so that $\qis(t')=\qis(t)$. (Possibly, $t'=t$ if $t\in\{t_1,t_2\}$.) For all integer $k$, we can find a point $t_k \in [t' - 1/k, t']$ such that $\Mi(t_k)$ is a leaf not identified with any other points by~$\sim$ (it suffices to select a point that is an increase point for neither functions of Lemma~\ref{lemip}). As a result, $\qis(t_k) \notin \qis([t_1,t_2])$ and $\disig(\qis(t_k),\qis(t))\to 0$ as $k\to\8$, so that $\qis(t) \in \partial (\qis([t_1,t_2]))$. The same argument shows that $ \Phi_{t_2}([w',\disig(s^\ooo,t_2)]) \subseteq \partial (\qis([t_1,t_2]))$ as well. To see that $\pii(\lhb \Mi(t_1), \Mi(t_2) \rhb ) \subseteq \partial (\qis([t_1,t_2]))$, we notice that
$$\pii\big(\lhb \Mi(t_1), \Mi(t_2) \rhb \big) = \qis\lp \Big\{ t \in [t_1,t_2] \, : \, C_\8^{e_t}\langle t\rangle = \inf_{r\in[t_1,t]} C_\8^{e_r}\langle r\rangle \text{ or } C_\8^{e_t}\langle t\rangle = \inf_{r\in[t,t_2]} C_\8^{e_r}\langle r\rangle \Big\} \rp$$
and apply a similar argument.

\smallskip

Conversely, let~$y \in \partial (\qis([t_1,t_2]))$ and $t\in [t_1,t_2]$ be such that $y=\qis(t)$ (note that $\qis([t_1,t_2])$ is a closed subset of~$\qis$ as the image of a closed subset of $[0,1]$ under the projection from $[0,1]$ into~$\qis$). We further suppose for the moment that~$y$ is not the root of~$\tau$ as this point needs a special argument. There exists a sequence $(r_k)_{k}$ of numbers in $[0,1]\bs [t_1,t_2]$ such that $\qis(r_k) \to \qis(t)$. Up to extraction, we may suppose that $r_k \to t' \in [0,1]\bs (t_1,t_2)$. Then $\qis(t')=\qis(t)$, so that either $t \eqt t'$, in which case $C_\8^{e_t}\langle t\rangle = \inf_{r\in[t_1,t]} C_\8^{e_r}\langle r\rangle$ or $C_\8^{e_t}\langle t\rangle = \inf_{r\in[t,t_2]} C_\8^{e_r}\langle r\rangle$ (we are necessarily in the case~\eqref{idsamefor} of the definition of~$\simeq$ from Section~\ref{secquot}, as we supposed that~$y$ is not the root of~$\tau$), or $t \not\eqt t'$ but $t\eqq t'$, in which case $\Lab(t) = \inf_{[t_1,t]} \Lab$ or $\Lab(t) = \inf_{[t,t_2]} \Lab$, by Proposition~\ref{propid}.

Finally, if the root of~$\tau$ belongs to $\partial( \qis([t_1,t_2]))$, it also belongs to the set we claim to be $\partial( \qis([t_1,t_2]))$ by a simple continuity argument. 
\end{pre}

We may now conclude that there is a unique geodesic from~$\bp$ to~$x$.

\begin{prop}\label{geodleaf}
The path~$\Phi_s$ is the only geodesic from~$\bp$ to~$x$.
\end{prop}

\begin{pre}
We consider a geodesic $\wp: [0,\disig(\bp,x)] \to \qis$ from~$\bp$ to~$x$. Let $\eps>0$ be such that $[s-\eps,s+\eps]$ is included in the interval coding~$\tau_a$. We take~$s_i$ and~$s_j$ satisfying the hypotheses of Lemma~\ref{passpas} with this value of~$\eps$ and we set $r\de \sup\{w\, : \, \wp(w) \notin \qis([s_i,s_j])\}$, so that $\wp(r)\in\partial(\qis([s_i,s_j]))$. Our choice of~$s_i$ and~$s_j$ ensures that $\wp(r) \notin \pii(\lhb a_i,a_j \rhb)$, hence, by Lemma~\ref{boundt1t2},
$$\wp(r) \in \Phi_{s_i}\big([w',\disig(s^\ooo,s_i)]\big) \cup \Phi_{s_j}\big([w',\disig(s^\ooo,s_j)]\big)\quad\text{ where }\quad w' \de \inf \big\{w\, : \, \Phi_{s_i}(w) \neq \Phi_{s_j}(w)\big\}.$$
As~$\Phi_{s_i}$ and~$\Phi_{s_j}$ are the only geodesics from~$\bp$ to~$x_i$ and~$x_j$, this yields that~$\wp$ restricted to $[0,r]$ is equal to either~$\Phi_{s_i}$ or~$\Phi_{s_j}$ restricted to the same interval. In particular, $\wp(w)=\Phi_{s_i}(w)=\Phi_{s_j}(w)$ for $w\in [0, w']$. See Figure~\ref{aaiajb2}. Conducting this reasoning with~$\Phi_s$, we obtain that $\wp(w)=\Phi_{s}(w)$ for $w\in [0, w']$. 

\begin{figure}[ht!]
		\psfrag{a}[][]{$x$}
		\psfrag{b}[][]{$b_\eta$}
		\psfrag{c}[][]{$c_\xi$}
		\psfrag{i}[][]{$a_i$}
		\psfrag{j}[][]{$a_j$}
		\psfrag{l}[][]{$l_\xi$}
		\psfrag{r}[r][r]{\textcolor{green!70!black}{$\wp(r)$}}
		\psfrag{u}[][]{$s$}
		\psfrag{s}[][]{$s_i$}
		\psfrag{t}[][]{$s_j$}
		\psfrag{e}[][]{$t_\eta$}
		\psfrag{p}[][]{\textcolor{green!70!black}{$\wp$}}
		\psfrag{1}[][]{\textcolor{blue}{$\Phi_{s_i}$}}
		\psfrag{2}[][]{\textcolor{blue}{$\Phi_{s_j}$}}
		\psfrag{x}[][]{$\xi$}
	\centering\includegraphics[width=55mm]{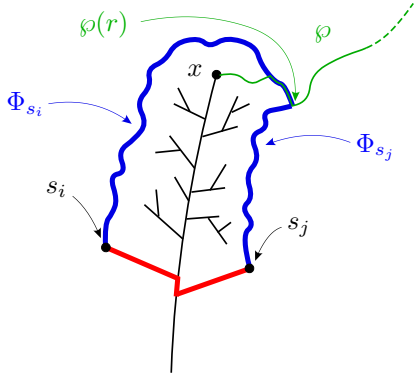}
	\caption{The geodesic~$\wp$ is ``trapped'' between~$\Phi_{s_i}$ and~$\Phi_{s_j}$.}
	\label{aaiajb2}
\end{figure}

As
$$\inf_{[s-\eps,s+\eps]}\Lab-\Lab(s^\ooo)\le w'=\inf_{[s_i\wedge s_j,s_i\vee s_j]}\Lab-\Lab(s^\ooo)\le \Lab(s)-\Lab(s^\ooo),$$
we conclude by continuity of~$\Lab$ that $\wp=\Phi_s$.
\end{pre}

\subsection{Geodesics from the base point to all points}

\begin{prop}\label{propallgeod}
Let $a\in\Mi$. The only geodesics from~$\bp$ to the point~$\pii(a)$ are the paths (all distinct)~$\Phi_s$, for all~$s$ such that $\Mi(s)=a$.
\end{prop}

\begin{pre}
We have already treated the case of leaves in Section~\ref{secgeodleaf}. The other cases are a little more intricate, but use exactly the same arguments, so that we only give a sketch of proof and leave the details to the reader.

If $a\in \Mi\bs \fl$ is a point of order~$2$, let $s < s'$ be such that $a=\Mi(s)=\Mi(s')$. For~$\eps$ small enough, the set $[s-\eps,s'+\eps]$ is included in the interval coding~$\tau_a$. For such an~$\eps$, following the ideas of Lemma~\ref{passpas}, we can find $s-\eps < s_i < s < s_j < s+\eps$ and $s'-\eps < s_{i'} < s' < s_{j'} < s'+\eps$ such that for all $b \in \lhb a_i,a_{j'} \rhb \cup \lhb a_{i'}, a_j \rhb$, we have $\disig(\bp,a) < \disig(\bp,b) + \disig(b,a)$.

We then follow the proof of Proposition~\ref{geodleaf}, replacing $\qis([s_i,s_j])$ with $\qis([s_i,s_j]\cup[s_{i'},s_{j'}])$, which is also a closed subset of~$\qis$, and whose boundary is the reunion of $\pii(\lhb a_i,a_{j'} \rhb \cup \lhb a_{i'}, a_j \rhb)$, the ranges of~$\Phi_{s_i}$ and~$\Phi_{s_j}$ from the point where they separate up to~$x_i$ and~$x_j$, as well as the ranges of~$\Phi_{s_{i'}}$ and~$\Phi_{s_{j'}}$ from the point where they separate up to~$x_{i'}$ and~$x_{j'}$. We eventually conclude that the only geodesic from~$\bp$ to~$\pii(a)$ are~$\Phi_s$ and~$\Phi_{s'}$ by the same technique.

All remaining cases are treated in a very similar fashion.
\end{pre}

We immediately deduce the following result by the same method as Le~Gall.

\begin{corol}[{\cite[Corollary~7.7]{legall08glp}}]\label{corconfgeod}
Almost surely, for every $\eps>0$, there exists $\eta\in\,(0,\eps)$ such that the following holds. Let~$\wp$ and~$\wp'$ be geodesics from~$\bp$ to points~$x$ and~$x'$ such that $\disig(\bp,x)\ge\eps$ and $\disig(\bp,x')\ge \eps$. Then $\wp(w)=\wp'(w)$ for every $w\in[0,\eta]$. 	
\end{corol}

\section{Remaining proofs}\label{secremproofs}

\subsection{Proof of Theorems~\ref{thmgeod} and~\ref{thmdimnc}}\label{secthm32}

\begin{pre}[Proof of Theorem~\ref{thmgeod}]
Note that, by~\eqref{boundary}, the indicator in Definition~\ref{deforder} is equal to $\un{\{\pii(a)\in\partial \qis\}}$. Assertion~($\ref{thmgeodi}$) is then a direct consequence of Proposition~\ref{propallgeod}. We now turn to the three remaining assertions. There is nothing to show in the case $(g,p)=(0,0)$ as all the sets in question are empty. If $(g,p)=(0,1)$, by the Jordan curve Theorem, all the cycles are homotopic to~$0$, so that we only need to see that $\pii(\cB)=\partial\qis$; this easily follows from~\eqref{boundary}. We thus suppose from now on that $(g,p)\notin\{(0,0),(0,1)\}$.

In this section, it will be convenient to ``remove'' the root edge dangling from~$\s_\8$. Recall that, as~$\s_\8$ is dominant, one extremity of the root~$e_*$ has degree~$1$ and the other one has degree~$3$. Let $\epsilon_1$, $\epsilon_2$, $\bar{\epsilon_2}$, $\epsilon_3$ be the sequence of successive half-edges in the contour order of~$f^\ooo$ such that $\{\epsilon_2, \bar{\epsilon_2}\}$ is the root edge of~$\s_\8$. Note that this implies that $\bar{\epsilon_1}$ is the half-edge following $\bar{\epsilon_3}$ in the contour of the incident face. We define the (nondominant,  unrooted) scheme~$\tilde\s_\8$ as the unrooted map obtained from~$\s_\8$ by removing the root edge $\{\epsilon_2, \bar{\epsilon_2}\}$ and merging $\{\epsilon_1, \bar{\epsilon_1}\}$ with $\{\epsilon_3, \bar{\epsilon_3}\}$. 
We associate with every half-edge $e\in\vec E(\s_\8)$ a half-edge~$\tilde e\in\vec E(\tilde\s_\8)$ as follows: if $e\in\{\epsilon_1,\epsilon_2,\bar{\epsilon_2},\epsilon_3\}$ then~$\tilde e$ is the merged edge, oriented from $\epsilon_1^-$ to $\epsilon_3^+$; if $e\in\{\bar{\epsilon_1},\bar{\epsilon_3}\}$ then~$\tilde e$ is the reverse of the previous half-edge and otherwise $\tilde e\de e$.

Recall from Section~\ref{secquot} that~$e_s$ is the half-edge of~$\s_\8$ ``visited'' at time~$s$. We also denote by $\bar\Phi_s$ the reverse of $\Phi_s$ and by $\wp\ooo\wp'$ the concatenation of~$\wp$ with~$\wp'$.

\begin{lem}\label{lemhomot}
Let $s\simeq s'$. The path $\Phi_s \ooo \bar\Phi_{s'}$ is homotopic to~$0$ if and only if $\tilde {e_s}= \tilde {e_{s'}}$.
\end{lem}

We postpone the proof of this lemma to the end of this section and conclude. Assertion~($\ref{thmgeodii}$) is shown as follows. First note that
$$\cB = \Mi\lp \lb s \, :\, C_\8^{e_s}\langle s \rangle = \underline { C_\8^{e_s}}\langle s \rangle \text{ and } e_s\notin\{e_*,\bar e_*\}\rb \rp.$$
This set splits into $\mathbf I\cup \mathbf B$ where
\begin{align}
\mathbf I&\de\Mi\lp \lb s \, :\, C_\8^{e_s}\langle s \rangle = \underline { C_\8^{e_s}}\langle s \rangle \text{ and } e_s\in\vec I(\s_\8)\bs\{e_*,\bar e_*\}\rb \rp,\label{setNC}\\
\mathbf B&\de\Mi\lp \lb s \, :\, C_\8^{e_s}\langle s \rangle = \underline { C_\8^{e_s}}\langle s \rangle \text{ and } e_s\in\vec B(\s_\8)\rb \rp.\notag
\end{align}
We already know by~\eqref{boundary} that $\pii(\mathbf B)=\partial\qis$ so that it will be sufficient to show that $\pii(\mathbf I)=\NC(\bp,\qis)$. If~$s$ lies in the set defining~$\mathbf I$, then there exists~$s'$ such that $s\simeq s'$ and $\tilde{e_{s'}}\neq\tilde{e_{s}}$ and as a consequence $\qis(s)\in \NC(\bp,\qis)$, by Lemma~\ref{lemhomot}. Conversely, if $\qis(s)\in \NC(\bp,\qis)$, then Lemma~\ref{lemhomot} implies the existence of~$s'$ such that $s\simeq s'$ and $\tilde {e_s}\neq \tilde {e_{s'}}$; as a result, the pair $\{s,s'\}$ satisfies~\eqref{idfacefor} or~\eqref{idnode} and thus $\Mi(s)\in \mathbf I$.

Assertion~($\ref{thmgeodiii}$) is also an easy consequence of Lemma~\ref{lemhomot}: a point belongs to the set in question if and only if it is of the form $\qis(s)=\qis(s')=\qis(s'')$ with $\tilde {e_s}$, $\tilde {e_{s'}}$, $\tilde {e_{s''}}$ pairwise distinct. These points are in one-to-one correspondence with the vertices of~$\tilde\s_\8$ that are not incident to any holes. Let $N(\s_\8)$ denotes their number. If $p=0$, these are all the vertices of~$\tilde\s_\8$: an easy application of Euler characteristic formula shows that $N(\s_\8)=4g-2$. If $p\ge 1$, the situation is a bit more complicated as~$N(\s_\8)$ actually depends on~$\s_\8$. It is a simple functional of~$\s_\8$ so that its distribution may be computed, as the distribution of~$\s_\8$ is given by Proposition~\ref{cvint}. Moreover, the distribution of~$\s_\8$ shows that~$\s_\8$ is equal to every dominant scheme of type $(p,1)$ with positive probability. This yields that the support of~$N(\s_\8)$ is merely the range of this functional. Let us now compute this range.

By Euler characteristic formula, $|V(\tilde\s_\8)|=4g+2p-2$. Because every hole of~$\s_\8$ is incident to at least one vertex, and because no vertices can be incident to more than one hole, we see that $N(\s_\8)\le 4g+p-2$. Now, for any given $k\in\{0,1,\ldots,4g+p-2\}$, it is not hard to construct a dominant scheme~$\s$ of type $(p,1)$ such that $N(\s)=k$ by using the three following operations (see Figure~\ref{domscheme}). Operation~(a) consists in appending an edge with an extremity incident to a degree $1$-hole: this increases the number~$N$ by one. Operation~(b) consists in adding a degree $2$-hole on some edge and operation~(c) consists in replacing a degree $\kappa$-vertex with a hole of degree~$\kappa$.

If $g\ge 1$ and $4g-2\le k\le 4g+p-2$, choose any dominant scheme of type $(0,1)$. At this stage, it has no holes and its number~$N$ is $4g-2$. Perform $k-(4g-2)$ times operation~(a) and $p-(k-(4g-2))$ times operation~(b) to obtain the desired result. If $g=0$, start with a dominant scheme of type $(1,1)$ and perform~$k$ times operation~(a) and $p-1-k$ times operation~(b). Finally, if $0\le k<4g-2$, take a (nondominant) scheme of type $(0,1)$ having exactly~$k+1$ degree $3$-vertices, one degree $1$-vertex and one degree $4g-k$-vertex. Perform operation~(c) on the degree $4g-k$-vertex and perform $p-1$ times operation~(b).

\begin{figure}[ht!]
		\psfrag{a}[][]{(a)}
		\psfrag{b}[][]{(b)}
		\psfrag{c}[][]{(c)}
	\centering\includegraphics[width=.95\linewidth]{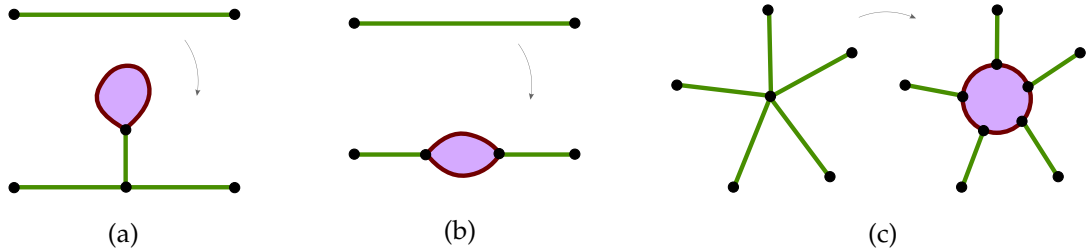}
	\caption{Adding a face in order to construct a dominant scheme~$\s$ for which $N(\s)=k$. \textbf{(a).} Adding one hole while increasing~$N$ by~$1$. \textbf{(b) and (c).} Adding one hole without changing~$N$.}
	\label{domscheme}
\end{figure}

Assertion~($\ref{thmgeodiv}$) directly follows from the observation that $\partial\qis \cap \NC(\bp,\qis)$ is in one-to-one correspondence with the set of vertices of~$\tilde\s_\8$ that are incident to a hole, so that its cardinality is equal to $|V(\tilde\s_\8)|-N(\s_\8)$.
\end{pre}

\begin{pre}[Proof of Theorem~\ref{thmdimnc}]
Here again, the Jordan curve Theorem allows us to directly conclude in the cases $(g,p)\in\{(0,0),(0,1)\}$, so that we exclude them. By the previous discussion, $\NC(\bp,\qis)=\pii(\mathbf I)$, where~$\mathbf I$ was defined by~\eqref{setNC}. In order to conclude that its Hausdorff dimension is a.s.\ equal to~$2$, we apply the method used in~\cite{bettinelli11slr} to compute the dimension of the boundary. In~\cite{bettinelli11slr}, the method was used to prove that the Hausdorff dimension of the boundary, that is, the image of the boundary edge of the scheme, was~$2$. In fact, the proof straightforwardly extends to any edge of the scheme and, as~$\mathbf I$ corresponds to a non-empty collection of scheme edges, the result follows.
\end{pre}

\begin{pre}[Proof of Lemma~\ref{lemhomot}]
Let $s\simeq s'$ and $a\de\Mi(s)=\Mi(s')$. If $s=s'$, the result is trivial so that we suppose $s\neq s'$. The condition $\tilde {e_s}= \tilde {e_{s'}}$ is equivalent to saying that one of the sets $\Mi(\st s {s'})$ or $\Mi(\st {s'} s)$ is a real tree. Indeed, if $\tilde {e_s}= \tilde {e_{s'}}$, then~$s$ and~$s'$ are instants corresponding to the same forest, so that $\Mi(\st s {s'})$ or $\Mi(\st {s'} s)$ is a real tree. Conversely, if $\tilde {e_s}\neq \tilde {e_{s'}}$ then~$a$ belongs to the backbone $\cB\subseteq\Mi$ and both $\Mi(\st s {s'})$ and $\Mi(\st {s'} s)$ contain a connected component of $\cB\bs\{a\}$. As a result, neither is a real tree.

Let us first suppose that $\tilde {e_s}= \tilde {e_{s'}}$ and, without loss of generality, let us suppose that $\Mi(\st s {s'})$ is a real tree. Let $w_0\de \inf\{w\, : \, \Phi_s(w)\neq\Phi_{s'}(w)\}$ and
\begin{equation}\label{rrr}
\rR\de \{\Phi_s(w),\Phi_{s'}(w)\,:  w_0\le w\le \disig(s^\ooo,s)\}.
\end{equation}
Plainly, $\Phi_s \ooo \bar\Phi_{s'}$ is homotopic to~$\rR$ so that it will be sufficient to show that~$\rR$ is homotopic to~$0$. We claim that~$\rR$ is equal to the boundary of $\big(\pii\big(\Mi(\st s {s'})\big),\disig\big)$: indeed, it is easy to check that if $\X\subseteq\X'$ are two surfaces with a boundary, then the boundary of~$\X$ is the union of its topological boundary with $(\partial \X')\cap\X$ and we obtain the claim thanks to Lemma~\ref{boundt1t2} (applied to $t_1\de s$ and $t_2\de s'$).

We will show that $\big(\pii\big(\Mi(\st s {s'})\big),\disig\big)$ is homeomorphic to a disk and the result will immediately follow. We can easily construct sequences~$(i_n)$ and~$(i'_n)$ such that $\m_n(i_n)=\m_n(i'_n)$, $i_n/2n\to s$, $i'_n/2n\to s'$ and the submap of~$\m_n$ composed of the elements visited between time~$i_n$ and time~$i'_n$ in the contour order is a tree~$\tr_n$. The simple discrete geodesic in~$\q_n$ from~$v_n^\ooo$ to~$i_n$ is the geodesic made of the arcs between the successive successors of~$i_n$ in the construction of Section~\ref{secrevmap}. The two simple geodesics toward~$i_n$ and~$i_n'$ may have a common part. Removing this common part, we are left with a simple loop, which we denote by~$\bm l_n$. The technique of Section~\ref{secconvms} implies that it is possible to isometrically embed the spaces of the pairs
$$\lp\Big( V(\q_{n_k}),\lp\frac 9 {8n_k}\rp^{1/4} d_{\q_{n_k}} \Big), \Big( V(\tr_{n_k}),\lp\frac 9 {8n_k}\rp^{1/4} d_{\q_{n_k}} \Big)\rp$$
as well as $\lp\big(\qis,\disig\big),\big(\pii\big(\Mi(\st s {s'})\big),\disig\big)\rp$ into the same metric space, in such a way that the first pair jointly converges toward the latter. It should be understood here that the metric is restricted to the considered sets. We cut the surface corresponding to~$\q_n$ introduced during Section~\ref{sectopo} along the loop $\bm l_n$. Clearly, the component corresponding to~$\tr_n$ is homeomorphic to a disk. Using \cite[Proposition~16]{bettinelli11slr}, we obtain that $\big(\pii\big(\Mi(\st s {s'})\big),\disig\big)$ is homeomorphic to a disk as well. (We need to check that the sequence~$(\bm l_n)_n$ is $0$-regular in the sense of \cite[Definition~9]{bettinelli11slr}; this is clear as~$\rR$ is a simple curve not reduced to a single point.)

\medskip

Conversely, let us suppose that neither $\Mi(\st s {s'})$ nor $\Mi(\st {s'} s)$ is a real tree. In particular, $a$ belongs to the backbone $\cB\subseteq\Mi$. Two cases are possible: either $\Mi(\st s {s'})\cap\Mi(\st {s'} s)=\{a\}$ or $\Mi(\st s {s'})\cap\Mi(\st {s'} s)\neq\{a\}$. Let us first suppose that we are in the latter case. Let us define~$\rR$ as previously by~\eqref{rrr}. Recall from Section~\ref{secsg} that the set $\pii^{-1}(\rR)$ is composed only of the point~$a$ and of leaves. We claim that this entails that the set $\qis\bs\rR$ is connected. Indeed, as $a\in\cB$, $\Mi\bs\pii^{-1}(\rR)$ may only have one or two connected components. If it has one component, the claim is immediate as $\qis\bs\rR$ is the image under~$\pii$ of $\Mi\bs\pii^{-1}(\rR)$. If it has two components, then~$a$ has to be of order~$3$ ($a$ corresponds to a node of~$\s_\8$ because $\Mi(\st s {s'})\cap\Mi(\st {s'} s)\neq\{a\}$). Let~$s''\notin\{s,s'\}$ be the third point such that~$\Mi(s'')=a$. Then the two images under~$\pii$ of the two connected components of $\Mi\bs\pii^{-1}(\rR)$ are connected and, for~$w$ large enough, $\Phi_{s''}(w)$ belongs to both images; the claim follows. Because~$\rR$ is not reduced to a single point, The Jordan curve Theorem entails that~$\rR$ cannot be homotopic to~$0$.

We now suppose that $\Mi(\st s {s'})\cap\Mi(\st {s'} s)=\{a\}$. The same method as above shows that both $\big(\pii\big(\Mi(\st s {s'})\big),\disig\big)$ and$\big(\pii\big(\Mi(\st {s'} s)\big),\disig\big)$ are surfaces with a boundary different from the disk (they contain at least a boundary component or a handle of~$\qis$), whose boundaries contain~$\rR$ as a component. As a result, $\rR$ cannot be homotopic to~$0$. 
\end{pre}

\subsection{Asymptotics for large quadrangulations}

The proof of Proposition~\ref{propdiscgeodi} is a straightforward adaptation of \cite[Section~9]{legall08glp}, so that we leave it to the reader and concentrate on the other proofs, which also use the same ideas but with a little extra work.

\begin{pre}[Proof of Proposition~\ref{prophmult}]
We argue by contradiction and suppose that there exists~$\eta>0$ such that, for infinitely many $n$'s, $\Pb(\exists v\in V(\q_n) : \HMult_{n}^{\eps_n}(v_n^\ooo,v)\ge 4) \ge\eta$. From the corresponding increasing sequence of integers, we extract a subsequence $(n_k)_{k\ge 0}$ along which the convergence of Theorem~\ref{cvthm} holds. We consider the simple loops made of edges of~$\q_{n_k}$ that are not homotopic to~$0$. Let $\delta_{n_k}$ denote the smallest diameter of these loops. An easy adaptation of \cite[Lemma~21]{bettinelli12tsl} shows that $\delta_{n_k}/{n_k}^{1/4}$ is uniformly bounded from below by an a.s.\ positive random variable~$\delta_0$.

We fix an integer~$n$ belonging to $(n_k)_{k\ge 0}$. We claim that for any~$v$ and~$v'$, $\HMult_{n}^\eps(v,v')\le \Mult_{n,\delta}^\eps(v,v')$ as soon as $\delta < \delta_n/(2n^{1/4})$. To see this, it is sufficient to see that every pair $\{\wp_n,\wp'_n\}$ of paths such that $\wp_n\ooo\bar\wp'_n$ is not homotopic to~$0$ satisfies $\dpath(\wp_n,\wp'_n)\ge \delta n^{1/4}$. We argue by contradiction and suppose that there exists a pair $\{\wp_n,\wp'_n\}$ satisfying the previous hypothesis but not the conclusion. It defines at least a simple loop~$\bm l_n$ that is not homotopic to~$0$. Let $v_1$, \dots, $v_j$ (resp.\ $v'_1$, \dots, $v'_j$) be the vertices of~$\bm l_n$ successively visited by~$\wp_n$ (resp.\ $\wp'_n$). For every $i \in\{1,\dots, j\}$, we consider a geodesic~$\gamma_i$ linking~$v_i$ to~$v'_i$. These geodesics have length smaller than $\delta n^{1/4}$ by hypothesis. For $i \in\{1,\dots, j-1\}$, the simple loop resulting of the concatenation of~$\gamma_i$, the edge visited by~$\wp'_n$ between time~$i$ and~$i+1$, $\bar\gamma_{i+1}$ and the edge visited by~$\bar \wp_n$ between time~$i+1$ and~$i$ has diameter smaller than $2\delta n^{1/4}$ so that it is homotopic to~$0$. As a result, the whole loop~$\bm l_n$ is homotopic to~$0$, which is a contradiction.

This entails that, for every $\delta >0$,
\begin{align*}
\Pb\big(\exists v\in V(\q_n) : \HMult_{n}^{\eps_n}(v_n^\ooo,v)\ge 4\big)\le \Pb\big(\delta_n\le 2\delta n^{1/4}\big) +\Pb\big(\exists v\in V(\q_n) : \Mult_{n,\delta}^{\eps_n}(v_n^\ooo,v)\ge 4\big)
\end{align*}
Using Proposition~\ref{propdiscgeodi}, we obtain the following contradiction:
\begin{equation*}
\eta \le\limsup_{k\to\8}\Pb\big(\exists v\in V(\q_{n_k}) : \HMult_{{n_k}}^{\eps_{n_k}}(v_{n_k}^\ooo,v)\ge 4\big)\le \inf_{\delta>0}\Pb\big(\delta_{0}\le 2\delta \big)=0.\qedhere
\end{equation*}
\end{pre}

Before we proceed, let us show our claim following Proposition~\ref{propdiscgeodii} about $\cB(\q_n,v_n^\ooo)\bs\partial V(\q_n)$. Let $v\in\cB(\q_n,v_n^\ooo)\bs\partial V(\q_n)$. We can find~$i$ and~$i'$ such that $v=\q_n(i)=\q_n(i')$ and the $i$-th and $i'$-th corners of~$\m_n$ do not belong to the same forest of the decomposition. It is a simple exercise to check that the simple discrete geodesics toward~$i$ and~$i'$ satisfy the required property, so that $\HMult_{n}^{0}(v_n^\ooo,v)\ge 2$.

\begin{pre}[Proof of Proposition~\ref{propdiscgeodii}]
We follow the lines of reasoning of~\cite[Section~9]{legall08glp}. We argue by contradiction and suppose that there exists~$\eta>0$ such that, for infinitely many $n$'s,
$$\Pb\big(\exists v\in V(\q_n)\, :\, d_{\q_n}(v,\cB(\q_n,v_n^\ooo)\bs\partial V(\q_n))\ge \delta n^{1/4},\, \HMult_{n}^{\eps_n}(v_n^\ooo,v)\ge 2\big) \ge\eta.$$
From the corresponding increasing sequence of integers, we extract a subsequence $(n_k)_{k\ge 0}$ along which the convergence of Theorem~\ref{cvthm} holds, and we let $(\qis, \disig)$ be the corresponding Brownian surface. We argue on the event
$$\limsup_{k\to\8} \big\{\exists v\in V(\q_{n_k})\, :\, d_{\q_{n_k}}(v,\cB(\q_{n_k},v_{n_k}^\ooo)\bs\partial V(\q_{n_k}))\ge \delta {n_k}^{1/4},\, \HMult_{{n_k}}^{\eps_{n_k}}(v_{n_k}^\ooo,v)\ge 2\big\},$$
whose probability is greater than~$\eta>0$. For infinitely many~$k$'s, we can find a vertex~$v_{n_k}\in V(\q_{n_k})$ that satisfies the above requirements; we let~$i_{n_k}$ be an integer such that $v_{n_k}=\q_{n_k}(i_{n_k})$. By compactness, we can find yet another subsequence along which $i_n/2n \to s$; we set $x\de\qis(s)$. Until further notice, we only consider integers~$n$ belonging to this subsequence. For every $\qis(s')\in\cB\bs\partial\qis$, we can easily construct a sequence $i_n'$ such that $\q_n(i'_n)\in \cB(\q_{n},v_{n}^\ooo)\bs\partial V(\q_{n})$ and $i'_n/2n \to s'$ along our subsequence. This implies that $d_{\q_n}(\q_n(i_n),\q_n(i'_n))\ge\delta n^{1/4}$, so that $\disig(x,\qis(s'))\ge (9/8)^{1/4}\delta$ and finally $\disig(x,\cB\bs\partial\qis)>0$. In particular, $x\notin \NC(\bp,\qis)$.

As $\HMult_{{n}}^{\eps_{n}}(v_{n}^\ooo,v_{n})\ge 2$, we can find two paths~$\wp_{n}$ and~$\wp'_{n}$ from~$v_{n}^\ooo$ to~$v_{n}$ with length less than $(1+\eps_{n})\,d_{\q_{n}}(v_{n}^\ooo,v_{n})$ and such that $\wp_{n}\ooo \bar{\wp}'_{n}$ is not homotopic to~$0$. There exist two integers~$w_n$ and $w_n'$ such that the part of~$\wp_n$ between times~$w_n$ and~$w_n'$ concatenated with the part of~$\bar\wp'_n$ between times~$w'_n$ and~$w_n$ make up a simple loop~$\bm l_n$ that is not homotopic to~$0$. We cut the surface of Section~\ref{sectopo} associated with~$\q_n$ along the loop~$\bm l_n$: we obtain either one or two surfaces with a boundary. As the number of possibilities for the topology of this or these surfaces is finite, there exist infinitely many $n$'s for which the topology is always the same. Restricting us to these~$n$'s, we see by the method of the previous section that this or these surfaces converge jointly with the surface associated with~$\q_n$ to a surface with the same topology and the boundary component corresponding to~$\bm l_n$ converges to a loop~$\bm l$ not homotopic to~$0$.

Now, let $0\le w\le \disig(\bp,x)$. We denote by~$\wp_n(w)$ (resp.\ $\wp'_n(w)$) the unique vertex of~$\wp_n$ (resp.\ of~$\wp'_n$) lying at distance $\lfloor(8n/9)^{1/4} w \wedge d_{\q_n}(v_n^\ooo,v_n)\rfloor$ from~$v_n^\ooo$. For each~$w$, we fix an extraction along which both~$\wp_n(w)$ and~$\wp_n'(w)$ converge and we denote by $\wp(w)$, $\wp'(w)\in\qis$ the limiting points. As~$\wp_n$ is an approximate geodesic, we see that $\disig(\wp(w),\wp(w'))=|w-w'|$, so that~$\wp$ is a geodesic from~$\bp$ to~$x$. Similarly, $\wp'$ is also a geodesic.

Up to further extraction, we may suppose that $w_n/(8n/9)^{1/4}\to w_\8$ and $w'_n/(8n/9)^{1/4}\to w'_\8$. As $w'_n-w_n$ is half the length of~$\bm l_n$, which, after scaling, does not converge to~$0$, we see that $w_\8<w'_\8$. On the one hand, for $w_\8<w <w'_\8$, we have $\wp_n(w)\in\bm l_n$ for~$n$ sufficiently large, so that $\wp(w)\in\bm l$. On the other hand, by definition, the vertex of~$\wp_n$ lying at distance~$w_n$ from~$v_n^\ooo$ is the same as the vertex of~$\wp'_n$ lying at distance~$w_n$ from~$v_n^\ooo$, so that $|\wp_n(w_\8)-\wp'_n(w_\8)|\to 0$ and thus $\wp(w_\8)=\wp'(w_\8)$. Moreover, as $\wp(w_\8)$ is a point of a geodesic that is not an end-point, it is a leaf by the discussion of Section~\ref{secsg}. As a result, $\wp$ and~$\wp'$ restricted to $[0,w_\8]$ are geodesics from~$\bp$ to a leaf, so that they coincide. Summing up, we constructed two geodesics~$\wp$ and~$\wp'$ from~$\bp$ to~$x$ such that $\wp(w)=\wp'(w)$ for $0\le w\le w_\8$ and $\wp([w_\8,w'_\8])\cup\wp'([w_\8,w'_\8]) \subseteq\bm l$. In order to obtain a contradiction with the fact that $x\notin \NC(\bp,\qis)$, we still need to see that $\wp\neq\wp'$. Note that this immediately entails that $w'_\8=\disig(\bp,x)$ as otherwise, $\wp(w'_\8)$ would be a leaf reachable from~$\bp$ by two distinct geodesics. To do this, let us take $y\in\bm l$. As~$\bm l_n$ converges toward~$\bm l$, $y$ is the limit of points of the form $\wp_n(r_n)$ or $\wp'_n(r_n)$ with $r_n\ge w_n/(8n/9)^{1/4}$. We first extract a subsequence along which $r_n\to r$. Considering the extraction along which both~$\wp_n(r)$ and~$\wp_n'(r)$ converge, we see that $y\in\{\wp(r),\wp'(r)\}$ so that actually $\wp([w_\8,\disig(\bp,x)])\cup\wp'([w_\8,\disig(\bp,x)]) =\bm l$, which concludes.
\end{pre}

\begin{pre}[Proof of Proposition~\ref{propdiscgeodiii}]
It is sufficient to show that
$$\lim_{\delta\to 0}\limsup_{n\to\8}\Pb\big(A_{n,\delta}^{\eps_n}(j)\big)\le \Pb(H\ge j)\le \lim_{\delta\to 0}\liminf_{n\to\8}\Pb\big(A_{n,\delta}^{\eps_n}(j)\big).$$
Let us focus on the first inequality. We argue by contradiction and suppose that there exist $\eta>0$ and $\delta>0$ such that, for infinitely many $n$'s, $\Pb(A_{n,\delta}^{\eps_n}(j))\ge \Pb(H\ge j)+\eta$. We extract from these $n$'s a subsequence $(n_k)_{k\ge 0}$ along which the convergence of Theorem~\ref{cvthm} holds, and we let $(\qis, \disig)$ be the corresponding Brownian surface. We argue on the event $\limsup_{k\to\8} A_{n_k,\delta}^{\eps_{n_k}}(j)$ whose probability is strictly greater than~$\Pb(H\ge j)$. Using the same technique as in the proof of Proposition~\ref{propdiscgeodii}, we can construct~$j$ points $x_1$, \ldots, $x_j\in\qis$ lying at distance greater than$(9/8)^{1/4}\delta$ one from another and belonging to the set appearing in~($\ref{thmgeodiii}$). As a result, we obtain that $\limsup_{k\to\8} A_{n_k,\delta}^{\eps_{n_k}}(j)\subseteq \{H\ge j\}$, which is a contradiction.

To show the second inequality, we do not extract any subsequences. We argue on the event $\{ H\ge j\}$, that is, the event that~$\tilde\s_\8$ has at least~$j$ vertices~$v_1$, \ldots, $v_j$ that are not incident to any holes (see Section~\ref{secthm32}). For~$n$ large enough, $\s_n=\s_\8$ so that these vertices correspond to vertices of~$V(\q_n)$. It is not hard to see that these vertices are such that $\HMult_{n}^{0}(v_n^\ooo,v_i)=3$. As a result, the event $A_{n,\delta}^{\eps_n}(j)$ holds, as long as $\min_{i\neq i'} d_{\q_n}(v_i,v_{i'})\ge \delta n^{1/4}$. The ``cactus bound'' (see~\cite{curienlgm13cactusI}) implies that the later condition holds for~$n$ large enough, provided that~$\delta$ is smaller than some a.s.\ positive random variable~$\Delta$, which is an explicit function of~$\s_\8$ and~$\Lab$. We just showed that
$$\Pb\big(\{\delta \le \Delta\} \cap \{H\ge j\}\big)\le\Pb\big(\liminf_{n\to\8}A_{n,\delta}^{\eps_n}(j)\big)$$
and we conclude thanks to Fatou's Lemma and by taking the limit $\delta\to 0$.
\end{pre}

The proof of Proposition~\ref{propdiscgeodiv} is very similar to the proof of Proposition~\ref{propdiscgeodiii}; we leave it to the reader.

\phantomsection
\addcontentsline{toc}{section}{Often used notation}
\label{secnot}
\printnomenclature[25mm]

\bibliographystyle{alpha}
\bibliography{main}
\end{document}